\documentclass[reqno,11pt]{amsart}
\usepackage[margin=1in]{geometry}
\usepackage{amsmath,amssymb,amsthm,cite}
\usepackage{enumitem}
\usepackage{mathrsfs}
\usepackage[colorlinks=true,allcolors=black]{hyperref}
\numberwithin{equation}{section}

\title[Steady Prandtl-Type Expansion with Strong Suction]{Validity of Steady Prandtl-Type Expansion\\ in the Half-Plane with Strong Suction}

\author[A. L. Mazzucato]{Anna L. Mazzucato}
\address{Department of Mathematics, Penn State University, University Park, PA 16802, USA.}
\email{alm24@psu.edu}

\author[D. Wang]{Dehua Wang}
\address{Department of Mathematics, University of Pittsburgh, Pittsburgh, PA 15260, USA.}
\email{dhwang@pitt.edu}

\author[W. Wang]{Weiqiang Wang}
\address{Department of Mathematics, University of Pittsburgh, Pittsburgh, PA 15260, USA.}
\email{wew179@pitt.edu}

\author[Y. Wang]{Yong Wang}
\address{Academy of Mathematics and Systems Science, Chinese Academy of Sciences, Beijing 100190, China; School of Mathematical Sciences, University of Chinese Academy of Sciences, Beijing 100049, China.}
\email{yongwang@amss.ac.cn}

\author[W. Wei]{Wei Wei}
\address{Department of Mathematics, University of Pittsburgh, Pittsburgh, PA 15260, USA.}
\email{wew144@pitt.edu}


\theoremstyle{plain}
\newtheorem{theorem}{Theorem}[section]
\newtheorem{lemma}{Lemma}[section]
\newtheorem{proposition}{Proposition}[section]
\newtheorem{corollary}{Corollary}[section]
\theoremstyle{definition}

\theoremstyle{remark}
\newtheorem{remark}{Remark}[section]

\begin{document}

\begin{abstract}
In this paper, a Prandtl-type expansion is justified and the vanishing viscosity limit is established for the two-dimensional steady incompressible Navier--Stokes equations in the half-plane with strong suction. Previous results on the Prandtl expansion with impermeable boundaries are valid for finite tangential length or infinite tangential length under additional assumptions on the Prandtl equations. Our result shows that strong suction prevents boundary-layer separation, allowing the expansion to hold globally in the tangential direction without any structural assumptions on the boundary-layer profiles. Two key ingredients of our proof are as follows: (1) A divergence-free multiplier involving a weight $y$ is constructed to obtain direct $L^2$ control of the remainder; (2) The singular terms involving the normal derivatives in the derivative estimates are controlled by Hardy-type inequalities derived by using the explicit exponential expression for the leading profile of the viscous boundary layer.
\end{abstract}

\maketitle

\section{Introduction}
In this paper, we study the vanishing viscosity limit of the steady Navier--Stokes equations in the half-plane $\Omega=\{(x,y)\in\mathbb{R}^2:y>0\}$ with suction boundary conditions:
\begin{equation}\label{nse}
\begin{cases}
u^\epsilon\partial_xu^\epsilon+v^\epsilon\partial_yu^\epsilon+\partial_xp^\epsilon=\epsilon\Delta u^\epsilon,\\
u^\epsilon\partial_xv^\epsilon+v^\epsilon\partial_yv^\epsilon+\partial_yp^\epsilon=\epsilon\Delta v^\epsilon,\\
\partial_xu^\epsilon+\partial_yv^\epsilon=0,\\
(u^\epsilon,v^\epsilon)(x,0)=(\beta(x),-\alpha(x)).
\end{cases}
\end{equation}
Here $(u^\epsilon,v^\epsilon)$, $p^\epsilon$, and $\epsilon>0$ denote the velocity, pressure, and viscosity, respectively. $(\beta(x),-\alpha(x))$ is the prescribed boundary data with $\alpha(x)>0$.

Formally, as $\epsilon\to 0$, system \eqref{nse} reduces to the steady Euler equations:
\begin{equation}\label{ee}
\begin{cases}
u^0\partial_xu^0+v^0\partial_yu^0+\partial_xp^0=0,\\
u^0\partial_xv^0+v^0\partial_yv^0+\partial_yp^0=0,\\
\partial_xu^0+\partial_yv^0=0,\\
v^0(x,0)=-\alpha(x).
\end{cases}
\end{equation}
Note that the tangential boundary condition is absent here, so in general there would be a tangential boundary mismatch
\begin{equation*}
\eta(x):=\beta(x)-u^0(x,0),
\end{equation*}
between the viscous and inviscid solutions, which gives rise to a viscous boundary layer.

For instance, suppose $(\beta,-\alpha)$ is constant. Then for any constant $\beta_0$,
\begin{equation*}
(u^{\epsilon},v^{\epsilon})=\bigl(\beta_0+(\beta-\beta_0)e^{-\frac{\alpha y}{\epsilon}},-\alpha\bigr)
\end{equation*}
is an exact solution of \eqref{nse}, which converges to the trivial solution of \eqref{ee} as $\epsilon\to 0$:
\begin{equation*}
(u^0,v^0)=(\beta_0,-\alpha).
\end{equation*}
Whenever $\beta_0\neq \beta$, the mismatch $\eta=\beta-\beta_0$ leads to the boundary-layer correction $\eta e^{-\alpha y/\epsilon}$, localized near the boundary on the scale $y=O(\epsilon)$.

The presence of boundary layers poses a major challenge in vanishing viscosity limit problems, and has long been one of the central topics in mathematical fluid mechanics. In 1904, Prandtl \cite{Prandtl-1904} proposed that incompressible flows in small viscosity regimes are governed by the Prandtl equations near the boundary, whereas they are approximated by the Euler equations away from the boundary, which is the Prandtl theory well-known to date. In the unsteady setting with the classical no-slip boundary condition, the Prandtl theory predicts the formation of a strong boundary layer, and the Prandtl expansions have been rigorously justified in the following settings: (i) analytic or Gevrey-class data \cite{GMM-2018,NN-2018,SC-1998a,SC-1998b,WW-2020,WWZ-2017,WWZ-2026}; (ii) initial data whose vorticity vanishes near the boundary \cite{FTZ-2018,Maekawa-2014}; and (iii) domains and initial data with suitable symmetries \cite{GKLMN-2019,LMT-2008,MT-2011}.
For characteristic boundaries with Navier slip, see \cite{BC-2023,BS-2012,CMR-1998,IP-2006,IS-2011,LNP-2005,MM-2017,NP-2011,Nguyen-2020,TWZ-2020,WXZ-2012,WYZ-2017,XX-2007,XY-1999,Zhong-2017}.
For noncharacteristic boundaries, see \cite{CC-2014,GHT-2012,GG-1998,GMWZ-2010,GK-2025,MM-2017,MWW-2026,TW-2002,YMMV-2026}.

In the steady setting with moving boundaries, Guo and Nguyen \cite{GN-2017} established the first result over a sufficiently short moving plate with a shear leading Euler flow. Iyer \cite{Iyer-2017,Iyer-2019-2020} subsequently treated a sufficiently narrow sector exterior to a rotating disk and a semi-infinite moving plate with a small velocity mismatch; see also \cite{Iyer-2019} for a non-shear leading Euler flow. For the no-slip boundary condition, Guo and Iyer \cite{GI-2021,GI-2023} justified the expansions over a sufficiently short stationary plate, and Gao and Zhang \cite{GZ-2023} extended this result to arbitrary finite plate lengths under monotonicity and concavity assumptions on the Prandtl profile. {Iyer and Masmoudi~\cite{IM-2026} justified the expansion on a semi-infinite plate for inflow data sufficiently close to a Blasius profile}; see also \cite{GZ-2024} for non-shear leading Euler flows over a sufficiently short plate, and \cite{GW-2025} for compressible flow. We also refer to \cite{GLMZ-2026} for three-dimensional steady axisymmetric Prandtl expansion.

If the domain is periodic in the tangential direction, the selection of the leading Euler flow is subtle. The Prandtl-Batchelor theory, which dates back to Prandtl \cite{Prandtl-1904} and was rediscovered by Batchelor \cite{Batchelor-1956}, predicts that a steady two-dimensional inviscid limit flow must have constant vorticity in an eddy. Gui, Xie, and Xu \cite{GXX-2026} recently removed the eddy assumption and proved the Prandtl-Batchelor theory for bounded connected domains with impermeable boundaries. Moreover, additional compatibility conditions such as the Batchelor-Wood criterion \cite{Wood-1957} are required for the leading Euler flow to guarantee the solvability of the Prandtl equations. Justification of the Prandtl-Batchelor theory through Prandtl expansions was made by Fei-Gao-Lin-Tao \cite{FGLT-2023,FGLT-2024} over a disk and an annulus, and Pan-Zhao \cite{PZ-2025} over the exterior of a disk.

When the domain is infinite in the tangential direction, potential boundary-layer separation happens, which causes difficulties in the validity of the Prandtl expansions. To avoid this scenario, most aforementioned results in the steady setting are established either over a short tangential extent or under certain structural assumptions on the Prandtl profile. In this work, we address this issue by imposing strong suction at the boundary. Indeed, experimentally and computationally, a sufficiently strong suction rate stabilizes the boundary layer \cite{INP-2007,MK-2000}. Also in engineering and aerodynamics, suction is a widely used technique to prevent boundary-layer separation \cite[Chapter 14]{SG-2017}.
{Recently, Li and Wang \cite{LW-2026} proved the asymptotic stability of steady Prandtl equations around the suction profile, $(u^\epsilon,v^\epsilon)=(1-e^{-\alpha y/\sqrt{\epsilon}},-\sqrt{\epsilon}\alpha)$ in our notation, where $\alpha>0$ is a constant.}
These results motivate the present study to justify the viscous boundary-layer expansion and establish the vanishing viscosity limit in the whole half-plane with strong suction.

To achieve this goal, we assume that, for some constants $\alpha_0>0$ and $s>9/2$,
\begin{equation}\label{reg}
\begin{gathered}
\alpha=\alpha_0+\tilde{\alpha},\quad\tilde{\alpha}\in H^s(\mathbb{R}),\quad \beta\in L^\infty(\mathbb{R}),\quad \beta'\in H^{s-1}(\mathbb{R}).
\end{gathered}
\end{equation}
We consider Euler flows dominated by suction, in the sense that there exists a constant $\lambda_0>0$ such that
\begin{equation}\label{alpha0}
-v^0\geq \lambda_0,\quad\lambda_0-\frac{1}{2}\sqrt{\lambda_0^2+(\|u^0\|_{L^\infty(\overline{\Omega})}+\|\eta\|_{L^\infty(\mathbb{R})})^2}>\|y\nabla(u^0,v^0)\|_{L^\infty(\Omega)}+\|\eta\|_{L^\infty(\mathbb{R})};
\end{equation}
We shall show that such Euler flow exists and provide an explicit characterization in Lemma \ref{0th-int-lem}. {From \eqref{alpha0} and the boundary condition $\eqref{ee}_4$, we have $\alpha(x)=-v^0(x,0)\geq \lambda_0$, which means the suction speed on the boundary is uniform over $x\in\mathbb{R}$.}

In this setting, our main result is as follows.
\begin{theorem}\label{main-thm}
Under condition \eqref{reg}, let $(u^0,v^0,p^0)$ be the Euler flow solving \eqref{ee} and satisfying \eqref{alpha0}. Then, the Navier-Stokes equations \eqref{nse} admit a strong solution $(u^\epsilon,v^\epsilon,p^\epsilon)$ for all sufficiently small $\epsilon>0$, satisfying
\begin{equation*}
\|(u^\epsilon-u^0-\eta e^{-\alpha y/\epsilon},v^\epsilon-v^0)\|_{C(\overline{\Omega})}\leq C\epsilon,
\end{equation*}
where $C>0$ is a constant independent of $\epsilon$.
\end{theorem}

\begin{remark}\leavevmode
\begin{enumerate}
\item In general, the uniqueness of the strong solution for the Navier-Stokes equations \eqref{nse} established in Theorem \ref{main-thm} is not known; see also \cite[Sec. XII.2, pp. 803--804]{Galdi-2011}.

\item
Since the $x$-direction ranges over the whole real line, condition \eqref{reg} is imposed to provide certain integrability of the approximate solution.
The inequality \eqref{alpha0} describes an Euler flow in which suction dominates the tangential velocity, boundary-layer mismatch, and the velocity variation.
For any fixed data $(\tilde{\alpha},\beta,\beta_0)$, condition \eqref{alpha0} holds for all sufficiently large background suction $\alpha_0$, even when $\|\tilde{\alpha}\|_{H^s}$ is large; see Lemma \ref{0th-int-lem}.

\item As aforementioned, strong suction can suppress boundary-layer separation, allowing us to establish Theorem \ref{main-thm} globally in the tangential direction.
\end{enumerate}
\end{remark}

We now outline the main ingredients of the proof.

(i) We begin by determining the interior profiles. Under \eqref{alpha0}, the leading Euler flow $(u^0,v^0)$ turns out to be irrotational. Indeed, the uniform suction $-v^0\geq \lambda_0$ ensures every backward characteristic reaches the far-field. Since the vorticity is transported along the characteristics and decays at infinity due to the bounded velocity variation $y\nabla(u^0,v^0)\in L^\infty(\Omega)$, it therefore vanishes identically. Irrotationality and incompressibility then imply that both velocity components are harmonic. This yields representations of $u^0$ and $v^0$ through the half-plane conjugate Poisson and Poisson kernels, respectively. The higher-order interior profiles satisfy the linearized Euler equations around the irrotational flow $(u^0,v^0)$. Under analogous boundedness assumptions, a similar argument shows these profiles are also irrotational and hence harmonic, yielding analogous representations; see Lemmas \ref{0th-int-lem}--\ref{1st-int-lem} for details.

(ii) Next, we construct the boundary-layer profiles. Since $v^0(x,0)=-\alpha(x)\neq0$, an asymptotic expansion in the scaled normal variable $z=y/\epsilon$ yields the following homogeneous second-order linear ODE for the leading boundary-layer profile $u_B^0$ (cf. \cite[Chapter 14]{SG-2017}):
\begin{equation*}
\begin{cases}
-\partial_{zz}u_B^0-\alpha(x)\partial_zu_B^0=0,\\
u_B^0(x,0)=\beta(x)-u^0(x,0)=\eta(x),
\end{cases}
\end{equation*}
which yields
\begin{equation*}
u_B^0(x,z)=\eta(x)e^{-\alpha(x)z}.
\end{equation*}
This profile is well-defined globally in the tangential direction, and decays exponentially in the normal direction due to $\alpha(x)\geq \lambda_0>0$.

Similarly, we solve inhomogeneous second-order linear ODEs with source terms determined by previously constructed profiles to obtain the higher-order tangential boundary-layer profiles. The incompressibility and normal momentum equations then yield the normal velocity and pressure boundary-layer profiles. These profiles explicitly combine regular x-dependent coefficients with factors decaying exponentially in $z$, a structure crucial for remainder estimates; see Section \ref{bdry} for details.

(iii) Let $(u^a,v^a,p^a)$ denote the approximate solution constructed by the viscous boundary-layer expansion in steps (i)-(ii). The remainder $(u,v,p):=(u^\epsilon,v^\epsilon,p^\epsilon)-(u^a,v^a,p^a)$ satisfies
\begin{equation}\label{remainder-eq0}
\begin{cases}
(u^a+u)\partial_xu+(v^a+v)\partial_yu+u\partial_xu^a+v\partial_yu^a+\partial_xp-\epsilon\Delta u=-R_u^a,\\
(u^a+u)\partial_xv+(v^a+v)\partial_yv+u\partial_xv^a+v\partial_yv^a+\partial_yp-\epsilon\Delta v=-R_v^a,\\
\partial_xu+\partial_yv=0,\\
(u,v)(x,0)=(0,0),\\
\lim\limits_{y\to\infty}(u,v)(x,y)=(0,0),
\end{cases}
\end{equation}
where $(R_u^a,R_v^a)$ is the residual. We truncate the domain to $\Omega_N:=(-N,N)\times(0,N)$ and consider the linearized problem:
\begin{equation}\label{trun-lin-eq0}
\begin{cases}
(u^a+\widehat{u}_N)\partial_xu_N+(v^a+\widehat{v}_N)\partial_yu_N+u_N\partial_xu^a+v_N\partial_yu^a+\partial_xp_N-\epsilon\Delta u_N=-R_u^a,\\
(u^a+\widehat{u}_N)\partial_xv_N+(v^a+\widehat{v}_N)\partial_yv_N+u_N\partial_xv^a+v_N\partial_yv^a+\partial_yp_N-\epsilon\Delta v_N=-R_v^a,\\
\partial_xu_N+\partial_yv_N=0,\\
(u_N,v_N)|_{\partial\Omega_N}=(0,0),
\end{cases}
\end{equation}
where $(\widehat{u}_N,\widehat{v}_N)$ is a given divergence-free vector field that vanishes on $\partial\Omega_N$ and is sufficiently small in $L^\infty(\Omega_N)$.

The energy estimates approach in \cite{GN-2017,FGLT-2023,FGLT-2024} does not apply directly here, since the leading normal Euler velocity $v_I^0$ does not vanish under the suction boundary condition. For instance, testing the momentum equations with $u_N$ and $v_N$ produces
\begin{equation*}
\Big|\int_{\Omega_N}v_N^2\partial_yv_I^0\,{\rm d}x{\rm d}y\Big| \leq\|\partial_yv_I^0\|_{L^\infty(\Omega_N)}\|v_N\|_{L^2(\Omega_N)}^2,
\end{equation*}
where the coefficient $\partial_yv_I^0$ need not vanish or be small. To control $\|v_N\|_{L^2(\Omega_N)}$, one may use Poincar\'e's inequality in $x$ or the Hardy inequality in $y$. However, Poincar\'e's inequality in $x$ introduces a constant depending on $N$ which cannot be uniformly controlled, while the Hardy inequality in $y$ introduces a coefficient $\|y^2\partial_yv_I^0\|_{L^\infty(\Omega_N)}$ which is not small and cannot be absorbed into the good terms coming from the viscous term.
We therefore need a different approach to obtain the $L^2$ estimate, instead of using the Hardy or Poincar\'e inequalities.

To this end, we introduce the stream function
\begin{equation*}
\psi_N=-\int_0^y u_N(x,s)\,{\rm d}s
\end{equation*}
and test the momentum equation with $(yu_N-\psi_N,yv_N)$. This divergence-free multiplier eliminates the pressure and transforms the advection term
\begin{equation*}
\Bigl((u^a+\widehat{u}_N)\partial_xu_N+(v^a+\widehat{v}_N)\partial_yu_N,\,
(u^a+\widehat{u}_N)\partial_xv_N+(v^a+\widehat{v}_N)\partial_yv_N\Bigr)
\end{equation*}
into a positive definite quadratic form in $(u_N,v_N)$:
\begin{align}
&\int_{\Omega_N}\Big[-\frac{3}{2}(v^a+\widehat{v}_N)u_N^2+(u^a+\widehat{u}_N)u_Nv_N-\frac{1}{2}(v^a+\widehat{v}_N)v_N^2\Big]\,{\rm d}x{\rm d}y \nonumber \\
&\geq\operatorname*{ess\,inf}_{\Omega_N}\Big(-(v^a+\widehat{v}_N)-\frac{1}{2}\sqrt{(v^a+\widehat{v}_N)^2+(u^a+\widehat{u}_N)^2}\Big)\|(u_N,v_N)\|_{L^2(\Omega_N)}^2. \nonumber
\end{align}

Meanwhile, the contribution of the advecting term involving the derivatives of the approximate solution $(u_N\partial_xu^a+v_N\partial_yu^a,\,u_N\partial_xv^a+v_N\partial_yv^a)$ is bounded by
\begin{equation*}
\left(\|y\nabla(u^0,v^0)\|_{L^\infty(\Omega_N)}+\|\eta\|_{L^\infty(\mathbb{R})}+C\epsilon\right)\|(u_N,v_N)\|_{L^2(\Omega_N)}^2.
\end{equation*}

Therefore, to close the $L^2$ estimate, we impose \eqref{alpha0}, which can be achieved by sufficiently large suction. Consequently, we obtain
\begin{equation*}
\|(u_N,v_N)\|_{L^2(\Omega_N)}^2+\epsilon\|\sqrt{y}\nabla(u_N,v_N)\|_{L^2(\Omega_N)}^2\leq C\|y(R_u^a,R_v^a)\|_{L^2(\Omega_N)}^2,
\end{equation*}
for some constant $C$ independent of $\epsilon$ and $N$.

Then for the $H^1$ estimate, we test the momentum equation with $(u_N,v_N)$. With the above $L^2$ estimate, the remaining difficulty is to control the singular normal derivative of the leading boundary-layer profile $u_B^0$. The explicit exponential form of $u_B^0$ is used to establish the following Hardy-type bound:
\begin{equation*}
\Big|\int_{\Omega_N}u_Nv_N\partial_y\bigl[u_B^0(x,y/\epsilon)\bigr]\,{\rm d}x{\rm d}y\Big|\leq\frac{\epsilon\|\eta\|_{L^\infty(\mathbb{R})}}{2\lambda_0}\|\partial_y(u_N,v_N)\|_{L^2(\Omega_N)}^2.
\end{equation*}
Since \eqref{alpha0} implies $\|\eta\|_{L^\infty(\mathbb{R})}<2\lambda_0$, this term is absorbed into the viscous dissipation. In conclusion, we obtain
\begin{align}
\epsilon\|\nabla(u_N,v_N)\|_{L^2(\Omega_N)}^2&\leq C\bigl(\|(u_N,v_N)\|_{L^2(\Omega_N)}^2+\|(R_u^a,R_v^a)\|_{L^2(\Omega_N)}\|(u_N,v_N)\|_{L^2(\Omega_N)}\bigr) \nonumber \\
&\leq C\bigl(\|y(R_u^a,R_v^a)\|_{L^2(\Omega_N)}^2+\|(R_u^a,R_v^a)\|_{L^2(\Omega_N)}\|y(R_u^a,R_v^a)\|_{L^2(\Omega_N)}\bigr), \nonumber
\end{align}
for some constant $C$ independent of $\epsilon$ and $N$; See Sections \ref{L^2} and \ref{H^1} for details.

(iv) To complete the proof, we establish the solvability of the remainder equations \eqref{remainder-eq0}. On each truncated domain, using uniform $H^1$ \textit{a priori} estimates, the well-posedness of the truncated linearized problem \eqref{trun-lin-eq0} follows from the standard Lax--Milgram theorem combined with G{\aa}rding inequality and the Fredholm alternative. Next, the uniform-in-$N$ Stokes estimate and Agmon inequality upgrade the regularity to $H^2$ and $L^\infty$. Then the existence of the truncated nonlinear problem follows from Schauder's fixed point theorem. Finally, the uniform-in-$N$ estimates and standard compactness arguments allow us to pass to the limit $N\to\infty$, obtaining a solution to the original global nonlinear problem \eqref{remainder-eq0} with the desired estimates; see Section \ref{exist} for details.

\medskip

The rest of the paper is organized as follows. In Section 2, we construct the approximate solution through the viscous boundary-layer expansion. The explicit formulas of the interior and boundary-layer profiles and their properties will be shown. In Section 3, we establish {\it a priori} uniform $H^2$ estimates for the solution of a linear remainder equation over the truncated domain $\Omega_N=(-N,N)\times (0,N)$. Then, the existence and uniform estimates for the solution of the nonlinear remainder equations and its limit $N\to \infty$ will be established through the Lax-Milgram theorem combined with the Fredholm alternative theorem and compactness arguments in Section 4. The proof of Theorem 1.1 will be given in Section 5.

\medskip

\textbf{Notation}: For a domain $D\subset\mathbb{R}^d$, $d=1,2$, $L^p(D)$ ($1\leq p\leq \infty$) denotes the standard Lebesgue space over $D$ with the norm $\|\cdot\|_{L^p(D)}$.
{For $s>0$, $H^s(D)$ denotes the standard $L^2$-based Sobolev space (see Triebel, \emph{Theory of Function Spaces II}~\cite{Triebel-1992}).}
$H_{0}^s(D)$ denotes the closure of all smooth compactly supported $C_c^\infty(D)$ in $H^s(D)$. $H_{\mathrm{loc}}^s(D)$ denotes the set of functions belonging to $H^s(D')$ for every relatively compact subdomain $D'$ of $D$. In this paper, $C(D)$ denotes the set of all bounded continuous functions over $D$ with the supremum norm, $C^{0,\gamma}_{b}(D)$ with $0<\gamma<1$ denotes the set of all bounded and H\"{o}lder continuous functions over $D$. $C$ denotes a generic constant independent of $\epsilon$, whose value may change from line to line.

\section{Approximate Solution}\label{approx}
In this section, we construct the approximate solution via a Prandtl-type expansion. Standard asymptotic analysis near the boundary (cf. \cite[Chapter 14]{SG-2017}) shows that the thickness of the viscous boundary layer under the suction boundary condition is $O(\epsilon)$. Thus, formally, one has the following expansion:
\begin{align}
u^\epsilon&\sim (u_I^0(x,y)+u^0_B(x,y/\epsilon))+\epsilon(u_I^1(x,y)+u^1_B(x,y/\epsilon))+\cdots, \nonumber \\
v^\epsilon&\sim (v_I^0(x,y)+v^0_B(x,y/\epsilon))+\epsilon(v_I^1(x,y)+v^1_B(x,y/\epsilon))+\cdots, \nonumber \\
p^\epsilon&\sim (p_I^0(x,y)+p^0_B(x,y/\epsilon))+\epsilon(p_I^1(x,y)+p^1_B(x,y/\epsilon))+\cdots. \label{approx-profile1}
\end{align}
Here the subscripts $I$ and $B$ denote the interior and boundary-layer profiles, respectively.

\subsection{Interior profiles}\label{int}
We first study the well-posedness of the interior profiles.

\subsubsection{Zeroth-order interior profile}\label{0th-int-sec}
Substituting the expansion \eqref{approx-profile1} into the Navier-Stokes equations \eqref{nse} and collecting terms of order $O(\epsilon^0)$, we find that the zeroth-order interior profile $(u_I^0,v_I^0,p_I^0)$ is precisely the Euler flow $(u^0,v^0,p^0)$ governed by \eqref{ee}. Under \eqref{alpha0}, the following lemma determines this profile explicitly up to some constants in $u^0$ and $p^0$.

\begin{lemma}\label{0th-int-lem}
Given any far-field state $\lim_{y\to\infty}u^0=\beta_0$, under \eqref{reg}, for sufficiently large $\alpha_0$ with $\tilde{\alpha}$ and $\beta$ fixed, the Euler equations \eqref{ee} admit a unique solution $(u^0,v^0,p^0)$ satisfying \eqref{alpha0}, up to an additive constant in $p^0$:
\begin{equation}\label{0th-int-sol}
u^0=\beta_0-\mathcal{Q}_y\ast\tilde{\alpha},\quad v^0=-\alpha_0-\mathcal{P}_y\ast\tilde{\alpha},\quad 
p^0=-\frac{1}{2}\Bigl((u^0)^2+(v^0)^2\Bigr).
\end{equation}
Here $\ast$ denotes convolution in $x$, and $\mathcal{P}_y(x)$ and $\mathcal{Q}_y(x)$ are the Poisson and conjugate Poisson kernels for the upper half-plane, respectively:
\begin{equation*}
\mathcal{P}_y(x)=\frac{1}{\pi}\frac{y}{x^2+y^2},\qquad\mathcal{Q}_y(x)=\frac{1}{\pi}\frac{x}{x^2+y^2}.
\end{equation*}
Moreover, $(u^0,v^0)$ satisfies the far-field limit
\begin{equation}\label{far-field-sol}
\lim\limits_{y\to\infty}(u^0,v^0)(x,y)=(\beta_0,-\alpha_0),
\end{equation}
and the estimates
\begin{equation}\label{regular}
\|(u^0-\beta_0,v^0+\alpha_0)\|_{L^\infty(\Omega)}+\|\nabla(u^0,v^0)\|_{H^{s-1/2}(\Omega)}+\|y\nabla(u^0,v^0)\|_{L^\infty(\Omega)}\leq C\|\tilde{\alpha}\|_{H^s(\mathbb{R})}.
\end{equation}
\end{lemma}

\begin{remark}\label{reg-eta}
Let $\mathcal{H}$ denote the Hilbert transform, which is defined through:
\begin{equation*}
\mathscr{F}(\mathcal{H}(f))(\xi)=-i\operatorname{sgn}(\xi)\mathscr{F}(f)(\xi),
\end{equation*}
where $\mathscr{F}$ denotes the Fourier transform. It follows from \cite[p.~324, (5.1.36)]{Grafakos-2014} that $\mathcal{Q}_y=\mathcal{H}(\mathcal{P}_y)$. Moreover, the multiplier representation shows that $\mathcal{H}$ is an isometry on $H^r(\mathbb{R})$ for every $r\in\mathbb{R}$. Thus, \eqref{0th-int-sol} gives the trace of $u^0(x,y)$ on $y=0$:
\begin{equation*}
u^0(x,0)=\beta_0-\mathcal{H}(\tilde{\alpha})(x).
\end{equation*}
The tangential boundary mismatch is therefore
\begin{equation*}
\eta(x)=\beta(x)-u^0(x,0)=\beta(x)-\beta_0+\mathcal{H}(\tilde{\alpha})(x).
\end{equation*}
By \eqref{reg} and Sobolev embedding $H^s(\mathbb{R})\hookrightarrow L^\infty(\mathbb{R})$, we have
\begin{equation*}
\eta\in L^\infty(\mathbb{R}),\quad\eta'\in H^{s-1}(\mathbb{R}).
\end{equation*}
\end{remark}

\begin{proof}
We proceed in several steps.

\noindent\textit{1}.
The condition $-v^0\geq \lambda_0$ in \ref{alpha0} indicates each backward characteristic reaches the far-field.

Indeed, let $(X(t),Y(t))$ be the backward characteristic starting from the initial point $(x_0,y_0)\in\Omega$:
\begin{equation*}
\begin{cases}
\frac{d}{dt}(X(t),Y(t))=-(u^0(X(t),Y(t)),v^0(X(t),Y(t))),\\
(X(0),Y(0))=(x_0,y_0).
\end{cases}
\end{equation*} 
It follows from the condition $y\nabla(u^0,v^0)\in L^\infty(\Omega)$ in \eqref{alpha0} that $\nabla(u^0,v^0)$ is bounded for $y\geq y_0>0$, so $(u^0,v^0)$ is Lipschitz there, and hence $(X(t),Y(t))$ is well-defined for $t\geq 0$.

Observing that
\begin{equation*}
\frac{d}{dt}Y(t)=-v^0(X(t),Y(t))\geq \lambda_0,
\end{equation*}
one has
\begin{equation*}
Y(t)\geq y_0+\lambda_0t.
\end{equation*}
In particular, this shows
\begin{equation}\label{limy}
\lim\limits_{t\to\infty}Y(t)=\infty.
\end{equation}

\noindent\textit{2}.
Let $\omega^0(x,y)=\partial_xv^0-\partial_yu^0$ be the vorticity. A direct calculation yields
\begin{equation*}
u^0\partial_x\omega^0+v^0\partial_y\omega^0=0,
\end{equation*}
which says that $\omega^0$ is transported along the characteristics:
\begin{equation}\label{vort}
\omega^0(X(t),Y(t))=\omega^0(x_0,y_0),\quad\forall t\geq 0.
\end{equation}

Since $y\nabla(u^0,v^0)$ is bounded, $\omega^0(x,y)\to 0$ as $y\to \infty$. Sending $t\to\infty$ in \eqref{vort} and using \eqref{limy}, we obtain $\omega^0(x,y)\equiv 0$ in $\Omega$, so the flow is irrotational.

\noindent\textit{3}.
We claim that the solution $(u^0,v^0)$ satisfying \eqref{alpha0}, if it exists, is unique up to a constant in $u^0$.

Indeed, irrotationality and incompressibility imply that $u^0$ and $v^0$ are harmonic. Suppose $(\tilde{u},\tilde{v})$ is the difference of two solutions satisfying \eqref{alpha0}. Then, $\tilde{u}$ and $\tilde{v}$ are harmonic, $\tilde{v}(x,0)=0$, and it follows from $u^0\in L^\infty(\Omega)$ in \eqref{alpha0} that $\tilde{u}$ is bounded. By irrotationality, $\tilde{u}$ satisfies
\begin{equation*}
\partial_y\tilde{u}(x,0)=\partial_x\tilde{v}(x,0)=0.
\end{equation*}
An even reflection of $\tilde{u}$ across $y=0$ gives a bounded harmonic function in the whole plane, which must be constant by Liouville's theorem. 
Then by irrotationality and incompressibility, $\tilde{v}$ satisfies
\begin{equation*}
\partial_x\tilde{v}=\partial_y\tilde{u}=0,\qquad\partial_y\tilde{v}=-\partial_x\tilde{u}=0.
\end{equation*}
As a result, $\tilde{v}$ is zero.

\noindent\textit{4}. 
In particular, $v^0$ satisfies
\begin{equation*}
\begin{cases}
\Delta v^0=0,\\
v^0(x,0)=-\alpha(x).
\end{cases}
\end{equation*}
A solution is given by the Poisson extension of the boundary data:
\begin{equation*}
v^0=-\mathcal{P}_y\ast\alpha=-\alpha_0-\mathcal{P}_y\ast\tilde{\alpha}=-\alpha_0-\frac{1}{\pi}\int_{\mathbb{R}}\frac{y}{(x-s)^2+y^2}\tilde{\alpha}(s)\,{\rm d}s.
\end{equation*}

It follows from incompressibility, irrotationality, and the identities
\begin{equation*}
\partial_x\mathcal{Q}_y=-\partial_y\mathcal{P}_y,\qquad\partial_y\mathcal{Q}_y=\partial_x\mathcal{P}_y,
\end{equation*}
that
\begin{equation*}
\partial_xu^0=-\partial_yv^0=\partial_y\mathcal{P}_y\ast\tilde{\alpha}=-\partial_x\mathcal{Q}_y\ast\tilde{\alpha},\qquad\partial_yu^0=\partial_xv^0=-\partial_x\mathcal{P}_y\ast\tilde{\alpha}=-\partial_y\mathcal{Q}_y\ast\tilde{\alpha},
\end{equation*}
and hence
\begin{equation*}
u^0(x,y)=b_0-\mathcal{Q}_y\ast\tilde{\alpha}=b_0-\frac{1}{\pi}\int_{\mathbb{R}}\frac{x-s}{(x-s)^2+y^2}\tilde{\alpha}(s)\,{\rm d}s,
\end{equation*}
for some constant $b_0$. In fact, we shall see in step \textit{6} that this solution is bounded, so the uniqueness in step \textit{3} applies.

Irrotationality and the Euler equations imply $\nabla(p^0+\frac12((u^0)^2+(v^0)^2))=0$, yielding the pressure formula up to a constant.

\noindent\textit{5}. 
Since
\begin{equation}\label{Young}
\|\mathcal{P}_y\|_{L^2(\mathbb{R})}=\|\mathcal{Q}_y\|_{L^2(\mathbb{R})}=\frac{1}{\sqrt{2\pi y}},
\end{equation}
Young's inequality gives
\begin{equation}\label{sqrty-L^infty}
\|\mathcal{P}_y\ast\tilde{\alpha}\|_{L^\infty(\mathbb{R})},\|\mathcal{Q}_y\ast\tilde{\alpha}\|_{L^\infty(\mathbb{R})}\leq\frac{\|\tilde{\alpha}\|_{L^2(\mathbb{R})}}{\sqrt{2\pi y}}\to0,\quad\text{as }y\to\infty.
\end{equation}
Consequently, we have the far-field limit:
\begin{equation*}
\lim\limits_{y\to\infty}(u^0,v^0)(x,y)=\lim\limits_{y\to\infty}(b_0-\mathcal{Q}_y\ast\tilde{\alpha},-\alpha_0-\mathcal{P}_y\ast\tilde{\alpha})(x,y)=(b_0,-\alpha_0).
\end{equation*}
The given far-field state $\lim_{y\to\infty}u^0=\beta_0$ implies $b_0=\beta_0$.

\noindent\textit{6}.
Since $\|\mathcal{P}_y\|_{L^1(\mathbb{R})}=1$, H\"older's inequality and Sobolev embedding yield
\begin{equation*}
\|\mathcal{P}_y\ast\tilde{\alpha}\|_{L^\infty(\mathbb{R})}\leq\|\tilde{\alpha}\|_{L^\infty(\mathbb{R})}\leq C\|\tilde{\alpha}\|_{H^s(\mathbb{R})}.
\end{equation*}
On the other hand, using properties of the Hilbert transform in Remark \ref{reg-eta}, we have
\begin{equation*}
\|\mathcal{Q}_y\ast\tilde{\alpha}\|_{L^\infty(\mathbb{R})}=\|\mathcal{H}(\mathcal{P}_y)\ast\tilde{\alpha}\|_{L^\infty(\mathbb{R})}=\|\mathcal{P}_y\ast\mathcal{H}(\tilde{\alpha})\|_{L^\infty(\mathbb{R})}\leq C\|\mathcal{H}(\tilde{\alpha})\|_{H^s(\mathbb{R})}=C\|\tilde{\alpha}\|_{H^s(\mathbb{R})}.
\end{equation*}
Combining these estimates, we obtain
\begin{equation*}
\|(u^0-\beta_0,v^0+\alpha_0)\|_{L^\infty(\Omega)}\leq C\|\tilde{\alpha}\|_{H^s(\mathbb{R})}.
\end{equation*}

\noindent\textit{7}.
To estimate $y\nabla(u^0,v^0)$, we begin by calculating
\begin{align}
\|y\partial_y\mathcal{Q}_y\|_{L^1(\mathbb{R})}=\|y\partial_x\mathcal{P}_y\|_{L^1(\mathbb{R})}&=\frac{y}{\pi}\int_{\mathbb{R}}\Big|\partial_x\Big(\frac{y}{x^2+y^2}\Big)\Big|\,{\rm d}x \nonumber \\
&=\frac{2}{\pi}\int_{\mathbb{R}}\frac{y^2|x|}{(x^2+y^2)^2}\,{\rm d}x=\frac{2}{\pi}\int_{\mathbb{R}}\frac{|t|}{(t^2+1)^2}\,{\rm d}t=\frac{2}{\pi}, \nonumber \\
\|y\partial_x\mathcal{Q}_y\|_{L^1(\mathbb{R})}=\|y\partial_y\mathcal{P}_y\|_{L^1(\mathbb{R})}
&=\frac{y}{\pi}\int_{\mathbb{R}}\Big|\partial_y\Big(\frac{y}{x^2+y^2}\Big)\Big|\,{\rm d}x \nonumber \\
&=\frac{y}{\pi}\int_{\mathbb{R}}\frac{|x^2-y^2|}{(x^2+y^2)^2}\,{\rm d}x=\frac{1}{\pi}\int_{\mathbb{R}}\frac{|t^2-1|}{(t^2+1)^2}\,{\rm d}t=\frac{2}{\pi}. \nonumber
\end{align}
Young's inequality and Sobolev embedding then yield
\begin{equation*}
\|y\nabla(u^0,v^0)\|_{L^\infty(\Omega)}\leq\frac{4}{\pi}\|\tilde{\alpha}\|_{L^\infty(\mathbb{R})}\leq C\|\tilde{\alpha}\|_{H^s(\mathbb{R})}.
\end{equation*}

\noindent\textit{8}.
The Fourier transform of the Poisson extension $\mathcal{P}_y\ast f$ in $x$ is given by
\begin{equation*}
\mathscr{F}(\mathcal{P}_y\ast f)(\xi)=e^{-y|\xi|}\mathscr{F}(f)(\xi).
\end{equation*}
Thus, each $x$-derivative of $\mathcal{P}_y\ast f$ contributes a factor of $|\xi|$ in the Fourier space, and each $y$-derivative contributes a $-|\xi|$ factor.
Using irrotationality and incompressibility, a direct calculation gives
\begin{align}
\|\nabla u^0\|_{H^{s-1/2}(\Omega)}^2=\|\nabla v^0\|_{H^{s-1/2}(\Omega)}^2&\leq C\int_0^\infty\int_\mathbb{R}\langle\xi\rangle^{2s-1}|\xi|^2e^{-2y|\xi|}|\mathscr{F}(\tilde{\alpha})(\xi)|^2\,{\rm d}\xi{\rm d}y \nonumber \\
&\leq C\int_\mathbb{R}\langle\xi\rangle^{2s-1}|\xi||\mathscr{F}(\tilde{\alpha})(\xi)|^2\,{\rm d}\xi \nonumber \\
&\leq C\int_\mathbb{R}\langle\xi\rangle^{2s}|\mathscr{F}(\tilde{\alpha})(\xi)|^2\,{\rm d}\xi= C\|\tilde{\alpha}\|_{H^s(\mathbb{R})}^2. \label{Hs}
\end{align}

\noindent\textit{9}.
Lastly, we show that \eqref{alpha0} holds for sufficiently large $\alpha_0$ with $\tilde{\alpha},\beta,\beta_0$ fixed. Indeed, with $\tilde{\alpha}$ fixed, $\mathcal{P}_y\ast\tilde{\alpha}$ is bounded, so the lower bound $\lambda_0$ for $-v^0$ can be made arbitrarily large by increasing $\alpha_0$. Then with $\tilde{\alpha},\beta,\beta_0$ fixed, $u^0,y\nabla (u^0,v^0),\eta$ are bounded uniformly in $\alpha_0$, so \eqref{alpha0} holds as long as $\alpha_0$ is chosen sufficiently large. This completes the proof of Lemma \ref{0th-int-lem}.
\end{proof}

\subsubsection{First-order interior profile}\label{1st-int-sec}
Fix the zeroth-order interior profile $(u_I^0,v_I^0,p_I^0)$. Collecting terms of order $O(\epsilon)$, we find that the first-order interior profile $(u_I^1,v_I^1,p_I^1)$ satisfies the Euler equations linearized about $(u_I^0,v_I^0)$: 
\begin{equation}\label{1st-int}
\begin{cases}
u_I^0\partial_xu_I^1+v_I^0\partial_yu_I^1+u_I^1\partial_xu_I^0+v_I^1\partial_yu_I^0+\partial_xp_I^1=0,\\
u_I^0\partial_xv_I^1+v_I^0\partial_yv_I^1+u_I^1\partial_xv_I^0+v_I^1\partial_yv_I^0+\partial_yp_I^1=0,\\
\partial_xu_I^1+\partial_yv_I^1=0,\\
v_I^1(x,0)=-v_B^1(x,0)=:-g(x),\\
\lim\limits_{y\to\infty}u_I^1(x,y)=0.
\end{cases}
\end{equation}
Here $\Delta u_I^0=\Delta v_I^0=0$ due to the harmonic nature of the zeroth-order interior profile. The boundary condition means the first-order interior profile cancels the normal error caused by the first-order boundary-layer profile on the boundary.

The following lemma determines the first-order interior profile explicitly up to a constant in $p_I^1$.

\begin{lemma}\label{1st-int-lem}
Let $g\in H^{s-1}(\mathbb{R})$. Then the linearized Euler equations \eqref{1st-int} admit a unique solution $(u_I^1,v_I^1,p_I^1)$ in the class
$(u_I^1,v_I^1),\,y\nabla(u_I^1,v_I^1)\in L^\infty(\Omega)$, up to a constant in $p_I^1$:
\begin{equation}\label{1st-int-sol}
u_I^1=-\mathcal{Q}_y\ast g,\quad v_I^1=-\mathcal{P}_y\ast g,\quad p_I^1=-u_I^0u_I^1-v_I^0v_I^1+C.
\end{equation}
Moreover, $(u_I^1,v_I^1)$ satisfies the far-field limit
\begin{equation}\label{far-field-sol1}
\lim\limits_{y\to\infty}(u_I^1,v_I^1)(x,y)=(0,0),
\end{equation}
and the estimates
\begin{align}
\|(u_I^1,v_I^1)\|_{L^\infty(\Omega)}+\|\nabla(u_I^1,v_I^1)\|_{H^{s-3/2}(\Omega)}+\|y\nabla(u_I^1,v_I^1)\|_{L^\infty(\Omega)}&\leq C\|g\|_{H^{s-1}(\mathbb{R})}; \nonumber \\
\|\sqrt{y}(u_I^1,v_I^1)\|_{L^\infty(\Omega)}+\|\sqrt{y}\nabla(u_I^1,v_I^1)\|_{L^2(\Omega)}&\leq C\|g\|_{L^2(\mathbb{R})}. \label{est-sqrty}
\end{align}
\end{lemma}

\begin{proof}
Taking the curl of \eqref{1st-int} and using the fact that the zeroth-order interior profile is irrotational yields the transport equation for the first-order interior vorticity $\omega_I^1=\partial_xv_I^1-\partial_yu_I^1$:
\begin{equation*}
u_I^0\partial_x\omega_I^1+v_I^0\partial_y\omega_I^1=0.
\end{equation*}
The same argument as in Lemma \ref{0th-int-lem} shows that $(u_I^1,v_I^1)$ is irrotational and satisfies all requirements of Lemma \ref{1st-int-lem} except \eqref{est-sqrty}.

The same argument used to derive \eqref{sqrty-L^infty} gives the first term in \eqref{est-sqrty}:
\begin{equation*}
\|\sqrt{y}(u_I^1,v_I^1)\|_{L^\infty(\Omega)}\leq C\|g\|_{L^2(\mathbb{R})}.
\end{equation*}
For the second term, Plancherel's identity yields
\begin{align}
\|\sqrt{y}\nabla(u_I^1,v_I^1)\|_{L^2(\Omega)}^2&\leq C\int_0^\infty\int_\mathbb{R}y|\xi|^2e^{-2y|\xi|}|\mathscr{F}(g)(\xi)|^2\,{\rm d}\xi{\rm d}y \nonumber \\
&=\frac{C}{4}\int_\mathbb{R}|\mathscr{F}(g)(\xi)|^2\,{\rm d}\xi\leq C\|g\|_{L^2(\mathbb{R})}^2. \nonumber
\end{align}
This completes the proof of Lemma \ref{1st-int-lem}.
\end{proof}

\subsection{Boundary-layer profiles}\label{bdry}
Next we establish the well-posedness of boundary-layer profiles.

To begin with, we introduce the scaled normal variable $z=y/\epsilon$, and then expand the interior profiles and their tangential derivatives near the boundary:
\begin{align}
u_I^0(x,\epsilon z)&=u_I^0(x,0)+\epsilon z\partial_yu_I^0(x,\epsilon\zeta_u^0), \nonumber \\
v_I^0(x,\epsilon z)&=v_I^0(x,0)+\epsilon z\partial_yv_I^0(x,0)+\frac{\epsilon^2z^2}{2}\partial_y^2v_I^0(x,\epsilon\zeta_v^0), \nonumber \\
v_I^1(x,\epsilon z)&=v_I^1(x,0)+\epsilon z\partial_yv_I^1(x,\epsilon\zeta_v^1), \nonumber \\
\partial_xu_I^0(x,\epsilon z)&=\partial_xu_I^0(x,0)+\epsilon z\partial_{xy}u_I^0(x,\epsilon\zeta_{u,x}^0), \nonumber \\
\partial_xv_I^0(x,\epsilon z)&=\partial_xv_I^0(x,0)+\epsilon z\partial_{xy}v_I^0(x,\epsilon\zeta_{v,x}^0), \label{int-exp}
\end{align}
for some $\zeta_u^0,\zeta_v^0,\zeta_v^1,\zeta_{u,x}^0,\zeta_{v,x}^0\in(0,z)$.

\subsubsection{Equations of $(v_B^0,p_B^0)$}
Substituting the expansions \eqref{approx-profile1} and \eqref{int-exp} into the Navier-Stokes equations $\eqref{nse}_{2,3}$ then collecting terms of order $O(\epsilon^{-1})$, we obtain the equations for $(v_B^0,p_B^0)$:
\begin{equation}\label{vb1}
\partial_zp_B^0=-v_I^0\partial_zv_B^0,\quad \partial_zv_B^0=0
\end{equation}

Since boundary-layer terms are expected to vanish away from the boundary, we impose the far-field decay conditions:
\begin{equation}\label{vb2}
\lim\limits_{z\to \infty}v_B^0=\lim\limits_{z\to \infty}p_B^0=0.
\end{equation}

Together, \eqref{vb1} and \eqref{vb2} imply
\begin{equation}\label{vb3}
v_B^0(x,z)=p_B^0(x,z)=0.
\end{equation}

\subsubsection{Equations of $(u_B^0,v_B^1,p_B^1)$}\label{0th-bdry-sec}
Collecting terms of order $O(\epsilon^{-1})$ from $\eqref{nse}_1$ and terms of order $O(\epsilon^0)$ from $\eqref{nse}_{2,3}$, we obtain the equations for $u_B^0$ and $(v_B^1,p_B^1)$:
\begin{equation}\label{0th-bdry}
\begin{cases}
-\alpha(x)\partial_zu_B^0=\partial_z^2u_B^0,\\
u_B^0(x,0)=\beta(x)-u_I^0(x,0)=\eta(x),\\
\lim\limits_{z\to\infty}u_B^0=0,
\end{cases}
\end{equation}
and
\begin{equation}\label{pb1}
\begin{cases}
\partial_zv_B^1=-\partial_xu_B^0,\\
\partial_zp_B^1=-u_B^0\partial_xv_I^0(x,0)-v_I^0(x,0)\partial_zv_B^1+\partial_{zz}v_B^1,\\
\lim\limits_{z\to \infty}v_B^1=\lim\limits_{z\to \infty}p_B^1=0.
\end{cases}
\end{equation}

A direct calculation yields
\begin{equation}\label{u_B^0}
u_B^0(x,z)=\eta(x)e^{-\alpha(x)z}, 
\end{equation}
and
\begin{equation}\label{v_B1}
v_B^1(x,z)=e^{-\alpha(x) z}\Big[\frac{\eta(x)'}{\alpha(x)}-\eta\alpha(x)'\Big(\frac{z}{\alpha(x)}+\frac{1}{\alpha(x)^2}\Big)\Big],\quad p_B^1(x,z)=-\frac{2\alpha'(x)\eta(x)}{\alpha(x)}e^{-\alpha(x)z}.
\end{equation}

By Remark \ref{reg-eta} and the regularity assumption \eqref{reg}, we have
\begin{equation*}
u_B^0(\cdot,z)\in L^\infty(\mathbb{R}),\quad\partial_xu_B^0(\cdot,z),v_B^1(\cdot,z)\in H^{s-1}(\mathbb{R}).
\end{equation*}
In particular, the trace of $v_B^1(x,z)$ on $z=0$ is given by
\begin{equation}\label{g}
g(x):=v_B^1(x,0)=\frac{\eta'(x)}{\alpha(x)}-\frac{\alpha'(x)}{\alpha(x)^2}\eta(x)=\Big(\frac{\eta(x)}{\alpha(x)}\Big)'\in H^{s-1}(\mathbb{R}).
\end{equation}

\subsubsection{Equations of $(u_B^1,v_B^2)$}\label{1st-bdry-sec}
Collecting all the terms of order $O(\epsilon^0)$ from $\eqref{nse}_1$ and terms of order $O(\epsilon^1)$ from $\eqref{nse}_3$, we obtain the equations for $u_B^1$ and $v_B^2$:
\begin{equation}\label{1st-bdry}
\begin{cases}
-\alpha(x)\partial_zu_B^1-\partial_z^2u_B^1=f(x,z),\\
u_B^1(x,0)=-u_I^1(x,0),\\
\lim\limits_{z\to\infty}u_B^1(x,z)=0,
\end{cases}
\end{equation}
and
\begin{equation}\label{1st-bdry2}
\begin{cases}
\partial_zv_B^2=-\partial_xu_B^1,\\
\lim\limits_{z\to \infty}v_B^2(x,z)=0,
\end{cases}
\end{equation}
where the source term $f(x,z)$ is given by
\begin{align}
f(x,z)&=-(u_I^0(x,0)+u_B^0)\partial_xu_B^0-u_B^0\partial_xu_I^0(x,0)-(-z\partial_xu_I^0(x,0)+v_I^1(x,0)+v_B^1)\partial_zu_B^0 \nonumber \\
&=-\Bigl[\bigl(\eta'(\beta_0-\mathcal{H}(\tilde{\alpha}))-\eta\mathcal{H}(\tilde{\alpha}')+\alpha\eta g\bigr)-z\eta\bigl(\alpha'(\beta_0-\mathcal{H}(\tilde{\alpha}))+\alpha\mathcal{H}(\tilde{\alpha}')\bigr)\Bigr]e^{-\alpha z}-\frac{\alpha'\eta^2}{\alpha}e^{-2\alpha z}. \nonumber
\end{align}

A direct calculation yields
\begin{align}
u_B^1(x,z)&=\Bigl[\mathcal{H}(g)+z\Bigl(\frac{2\eta}{\alpha}\mathcal{H}(\tilde{\alpha}')-\beta g\Bigr)+\frac{z^2\eta}{2\alpha}\bigl(\alpha'(\beta_0-\mathcal{H}(\tilde{\alpha}))+\alpha\mathcal{H}(\tilde{\alpha}')\bigr)\Bigr]e^{-\alpha z} \nonumber \\
&+\frac{\alpha'\eta^2}{2\alpha^3}\bigl(e^{-2\alpha z}-e^{-\alpha z}\bigr), \label{u_B^1}
\end{align}
and
\begin{align}
v_B^2(x,z)&=\partial_x\Bigl\{\frac{e^{-\alpha z}}{\alpha}\Bigl[\mathcal{H}(g)+\Bigl(z+\frac{1}{\alpha}\Bigr)\Bigl(\frac{2\eta}{\alpha}\mathcal{H}(\tilde{\alpha}')-\beta g\Bigr) \nonumber \\
&\quad+\frac{\eta}{2\alpha}\Bigl(z^2+\frac{2z}{\alpha}+\frac{2}{\alpha^2}\Bigr)\Bigl(\alpha'(\beta_0-\mathcal{H}(\tilde{\alpha}))+\alpha\mathcal{H}(\tilde{\alpha}')\Bigr)\Bigr]+\frac{\alpha'\eta^2}{4\alpha^4}\bigl(e^{-2\alpha z}-2e^{-\alpha z}\bigr)\Bigr\}. \label{v_B^2}
\end{align}

By the regularity assumption \eqref{reg} and properties of the Hilbert transform, we have $u_B^1(\cdot,z)\in H^{s-1}(\mathbb{R})$ and $v_B^2(\cdot,z)\in H^{s-2}(\mathbb{R})$. In particular, the trace of $v_B^2(x,z)$ at $z=0$ is given by
\begin{equation*}
h(x):=v_B^2(x,0)=\partial_x\Bigl[\frac{\mathcal{H}(g)}{\alpha}-\frac{\beta g}{\alpha^2}+\frac{3\eta}{\alpha^3}\mathcal{H}(\tilde{\alpha}')+\frac{\alpha'\eta}{\alpha^4}\Bigl(\beta-\frac{5\eta}{4}\Bigr)\Bigr]\in H^{s-2}(\mathbb{R}).
\end{equation*}

\subsection{Construction of approximate solution}\label{corrector-sec}
Combining the interior and boundary-layer profiles, we construct the approximate solution:
\begin{align}
u^a(x,y)&=(u^0_I(x,y)+u^0_B(x,y/\epsilon))+\epsilon(u^1_I(x,y)+u^1_B(x,y/\epsilon))+\epsilon^2u^c(x,y), \nonumber \\
v^a(x,y)&=v^0_I(x,y)+\epsilon(v^1_I(x,y)+v^1_B(x,y/\epsilon))+\epsilon^2(v_B^2(x,y/\epsilon)+v^c(x,y)), \nonumber \\
p^a(x,y)&=p^0_I(x,y)+\epsilon(p^1_I(x,y)+p^1_B(x,y/\epsilon)). \label{approx-profile}
\end{align}
Here $(u^c,v^c)$ is an artificial corrector that forces $(u^a,v^a)$ to satisfy incompressibility and the boundary condition:
\begin{equation*}
\partial_xu^a+\partial_yv^a=0,\quad (u^a,v^a)(x,0)=(\beta(x),-\alpha(x)).
\end{equation*}
We therefore seek a corrector $(u^c,v^c)$ such that
\begin{equation}\label{u_c}
\partial_xu^c+\partial_yv^c=0,\quad(u^c,v^c)(x,0)=(0,-v_B^2(x,0))=(0,-h(x)).
\end{equation}

Let $\rho\in C_c^\infty([0,\infty))$ be a cutoff with $\rho(y)=1$ for $0\leq y\leq1$. We construct $(u^c,v^c)$ as follows:
\begin{equation}\label{corrector}
u^c(x,y)=\rho'(y)\int_0^\infty u_B^1(x,z)\,{\rm d}z,\quad v^c(x,y)=-\rho(y)h(x).
\end{equation}
It is straightforward to verify the construction satisfies \eqref{u_c}. Since $u_B^1(\cdot,z)\in H^{s-1}(\mathbb{R})$ and $h(x)\in H^{s-2}(\mathbb{R})$, we have $u^c(\cdot,y)\in H^{s-1}(\mathbb{R})$ and $v^c(\cdot,y)\in H^{s-2}(\mathbb{R})$.

Consequently, we arrive at the following proposition on the approximate solution.
\begin{proposition}\label{approx-profile-prop}
The approximate solution $(u^a,v^a,p^a)$ constructed by \eqref{approx-profile} is divergence-free and satisfies the boundary condition
\begin{equation*}
(u^a,v^a)(x,0)=(\beta(x),-\alpha(x)).
\end{equation*}
Moreover, for $i\in\{0,1\}$, $j\in\{1,2\}$, $m\in\{1,2,3\}$, $n\in\{1,2\}$, $l\in\{0,1,2\}$, and $k\in\mathbb{N}_0$, it holds that
\begin{align}
&(u_I^i,v_I^i)\in L^\infty(\Omega),\quad\nabla^m(u_I^0,v_I^0), \nabla^n(u_I^1,v_I^1)\in L^2(\Omega)\cap L^\infty(\Omega),\quad y\nabla(u_I^i,v_I^i)\in L^\infty(\Omega), \nonumber \\
&\sqrt{y}(u_I^1,v_I^1)\in L^\infty(\Omega),\quad\sqrt{y}\nabla(u_I^1,v_I^1)\in L^2(\Omega); \label{interior} \\
&z^k\nabla_{x,z}^l u_B^0\in L^\infty(\Omega),\quad z^k\nabla_{x,z}^l(\partial_xu_B^0,u_B^1,v_B^j)\in L^2(\Omega)\cap L^\infty(\Omega); \label{boundary} \\
&y^k\nabla_{x,y}^l(u^c,v^c)\in L^2(\Omega)\cap L^\infty(\Omega). \label{corr}
\end{align}
\end{proposition}

\begin{proof}
Estimates \eqref{interior} for the interior profiles follow directly from Lemmas \ref{0th-int-lem} and \ref{1st-int-lem} with Sobolev embedding $H^r(\Omega)\hookrightarrow L^\infty(\Omega)$ for $r>1$. Estimates \eqref{boundary} for the boundary-layer profiles follow from the explicit expressions \eqref{u_B^0}, \eqref{v_B1}, \eqref{u_B^1}, \eqref{v_B^2}, and the Sobolev embedding $H^r(\mathbb{R})\hookrightarrow L^\infty(\mathbb{R})$ for $r>1/2$. The same argument applied to \eqref{corrector} yields \eqref{corr}. Therefore, the proof of Proposition \ref{approx-profile-prop} is complete.
\end{proof}

\subsection{Residual}
With the approximate profile $(u^a,v^a,p^a)$ in hand, we finally estimate its residual:
\begin{equation*}
(R_u^a,R_v^a)=(u^a\partial_xu^a+v^a\partial_yu^a+\partial_xp^a-\epsilon\Delta u^a,u^a\partial_xv^a+v^a\partial_yv^a+\partial_yp^a-\epsilon\Delta v^a).
\end{equation*}

Using \eqref{ee}, \eqref{1st-int}, \eqref{int-exp}, \eqref{vb3}, \eqref{0th-bdry}, \eqref{pb1}, \eqref{1st-bdry}, and \eqref{1st-bdry2}, a careful computation yields
\begin{align}
R_u^a&=\epsilon\bigl[z\partial_yu_I^0(x,\epsilon\zeta_u^0)\partial_xu_B^0+zu_B^0\partial_{xy}u_I^0(x,\epsilon\zeta_{u,x}^0) \nonumber \\
&\quad+z\partial_yv_I^0(x,0)\partial_zu_B^1+\frac{z^2}{2}\partial_y^2v_I^0(x,\epsilon\zeta_v^0)\partial_zu_B^0+z\partial_yv_I^1(x,\epsilon\zeta_v^1)\partial_zu_B^0 \nonumber \\
&\quad+u_B^0\partial_xu_I^1+(u_I^0+u_B^0)\partial_xu_B^1+u_B^1\partial_xu_I^0+(u_I^1+u_B^1)\partial_xu_B^0 \nonumber \\
&\quad+v_B^1\partial_yu_I^0+(v_I^1+v_B^1)\partial_zu_B^1+(v_B^2+v^c)\partial_zu_B^0+\partial_xp_B^1-\partial_x^2u_B^0\bigr] \nonumber \\
&+\epsilon^2\bigl[\frac{z^2}{2}\partial_y^2v_I^0(x,\epsilon\zeta_v^0)\partial_zu_B^1 \nonumber \\
&\quad+(u_I^0+u_B^0)\partial_xu^c+(u_I^1+u_B^1)\partial_x(u_I^1+u_B^1)+u^c\partial_x(u_I^0+u_B^0) \nonumber \\
&\quad+v_I^0\partial_yu^c+(v_I^1+v_B^1)\partial_yu_I^1+(v_B^2+v^c)(\partial_yu_I^0+\partial_zu_B^1)-\partial_x^2u_B^1\bigr] \nonumber \\
&+\epsilon^3\Bigl[(u_I^1+u_B^1)\partial_xu^c+u^c\partial_x(u_I^1+u_B^1)+(v_I^1+v_B^1)\partial_yu^c+(v_B^2+v^c)\partial_yu_I^1-\Delta u^c\bigr] \nonumber \\
&+\epsilon^4\bigl[u^c\partial_xu^c+(v_B^2+v^c)\partial_yu^c\bigr], \nonumber
\end{align}
and
\begin{align}
R_v^a&=\epsilon\bigl[z\partial_yv_I^0(x,0)\partial_zv_B^1+zu_B^0\partial_{xy}v_I^0(x,\epsilon\zeta_{v,x}^0)+u_B^0\partial_xv_I^1+(u_I^0+u_B^0)\partial_xv_B^1+u_B^1\partial_xv_I^0 \nonumber \\
&\quad+v_I^0\partial_zv_B^2+v_B^1\partial_yv_I^0+(v_I^1+v_B^1)\partial_zv_B^1-\partial_z^2v_B^2\bigr] \nonumber \\
&+\epsilon^2\bigl[\frac{z^2}{2}\partial_y^2v_I^0(x,\epsilon\zeta_v^0)\partial_zv_B^1 \nonumber \\
&\quad+(u_I^0+u_B^0)\partial_x(v_B^2+v^c)+(u_I^1+u_B^1)\partial_x(v_I^1+v_B^1)+u^c\partial_xv_I^0 \nonumber \\
&\quad+v_I^0\partial_yv^c+(v_I^1+v_B^1)(\partial_yv_I^1+\partial_zv_B^2)+(v_B^2+v^c)(\partial_yv_I^0+\partial_zv_B^1)-\partial_x^2v_B^1\bigr] \nonumber \\
&+\epsilon^3\bigl[(u_I^1+u_B^1)\partial_x(v_B^2+v^c)+u^c\partial_x(v_I^1+v_B^1) \nonumber \\
&\quad+(v_I^1+v_B^1)\partial_yv^c+(v_B^2+v^c)(\partial_yv_I^1+\partial_zv_B^2)-\partial_x^2v_B^2-\Delta v^c\bigr] \nonumber \\
&+\epsilon^4\bigl[u^c\partial_x(v_B^2+v^c)+(v_B^2+v^c)\partial_yv^c\bigr]. \nonumber
\end{align}

Observe that the highest-order tangential derivatives of the boundary data appearing in the above residual expressions are of fourth order, which occur in $\partial_x^2v_B^2$ and $\partial_x^2v^c$. 
By Proposition \ref{approx-profile-prop} and a careful computation, we establish the following estimates for the residual:
\begin{equation}\label{res-est}
\|(R^a_u,R^a_v)\|_{L^2(\Omega)}\leq C\epsilon^{3/2},\qquad\|y(R^a_u,R^a_v)\|_{L^2(\Omega)}\leq C\epsilon^2,
\end{equation}
for some constant $C>0$ independent of $\epsilon$.

Indeed, for terms involving $u_B^0=\eta e^{-\alpha y/\epsilon}$, althrough $u_B^0\notin L^2(\Omega)$ since $\eta\notin L^2(\mathbb{R})$, they are a product of $u_B^0\in L^\infty(\Omega)$ and another function in $L^2(\Omega)$, and therefore lie in $L^2(\Omega)$. For the weighted estimate, the quadratic term $y(u_I^1,v_I^1)\cdot\nabla(u_I^1,v_I^1)$ requires particular attention, for which we have used \eqref{est-sqrty} to obtain
\begin{equation*}
\|\epsilon^2y(u_I^1,v_I^1)\cdot\nabla(u_I^1,v_I^1)\|_{L^2(\Omega)}\leq\epsilon^2\|\sqrt{y}(u_I^1,v_I^1)\|_{L^\infty(\Omega)}\|\sqrt{y}\nabla(u_I^1,v_I^1)\|_{L^2(\Omega)}\leq C\epsilon^2.
\end{equation*}
The other terms are estimated directly using Hölder's inequality and Proposition \ref{approx-profile-prop}.

\section{{\it A Priori} Estimates for Truncated Linear Problem}\label{est}
With the approximate profile $(u^a,v^a,p^a)$ constructed in \eqref{approx-profile}, we now turn to the {\it a priori} estimates for the remainder
\begin{equation*}
(u,v,p):=(u^\epsilon-u^a,v^\epsilon-v^a,p^\epsilon-p^a)
\end{equation*}
which satisfies the following equations:
\begin{equation}\label{remainder-eq}
\begin{cases}
(u^a+u)\partial_xu+(v^a+v)\partial_yu+u\partial_xu^a+v\partial_yu^a+\partial_xp-\epsilon\Delta u=-R_u^a,\\
(u^a+u)\partial_xv+(v^a+v)\partial_yv+u\partial_xv^a+v\partial_yv^a+\partial_yp-\epsilon\Delta v=-R_v^a,\\
\partial_xu+\partial_yv=0,\\
(u,v)(x,0)=(0,0),\\
\lim\limits_{y\to\infty}(u,v)(x,y)=(0,0).
\end{cases}
\end{equation}

To begin, we introduce the truncated domain:
\begin{equation*}
\Omega_N:=(-N,N)\times(0,N),\quad \text{for $N\in \mathbb{N}$},
\end{equation*}
and define the associated function spaces
\begin{align}
&\mathcal{V}_N:=\{(\varphi_u,\varphi_v)\in H_0^1(\Omega_N)\,:\,\partial_x\varphi_u+\partial_y\varphi_v=0\}, \nonumber \\
&\mathcal{X}_N:=\{(\varphi_u,\varphi_v)\in C(\overline{\Omega_N})\,:\,\partial_x\varphi_u+\partial_y\varphi_v=0,\varphi_u=\varphi_v=0\text{ on }\partial\Omega_N\}, \nonumber \\
&\mathcal{B}_N^\delta:=\{(\varphi_u,\varphi_v)\in \mathcal{X}_N\,:\,\|(\varphi_u,\varphi_v)\|_{C(\overline{\Omega_N})}\leq \delta\}. \nonumber
\end{align}
Here $\partial_x\varphi_u+\partial_y\varphi_v=0$ in the definition of $\mathcal{X}_N$ is understood in the distributional sense. Equipped with the sup-norm, $\mathcal{X}_N$ is a Banach space, and $\mathcal{B}_N^\delta$ is a closed convex subset of $\mathcal{X}_N$ for any $\delta>0$.

To justify integration by parts without {\it a priori} decay at infinity in the estimates, we first truncate the unbounded domain $\Omega$ to $\Omega_N$ and impose homogeneous Dirichlet boundary conditions, obtaining
\begin{equation}\label{trun-nonlin-eq}
\begin{cases}
(u^a+u_N)\partial_xu_N+(v^a+v_N)\partial_yu_N+u_N\partial_xu^a+v_N\partial_yu^a+\partial_xp_N-\epsilon\Delta u_N=-R_u^a,\\
(u^a+u_N)\partial_xv_N+(v^a+v_N)\partial_yv_N+u_N\partial_xv^a+v_N\partial_yv^a+\partial_yp_N-\epsilon\Delta v_N=-R_v^a,\\
\partial_xu_N+\partial_yv_N=0,\\
(u_N,v_N)=(0,0)\text{ on }\partial\Omega_N.
\end{cases}
\end{equation}
Then, given $(\widehat{u}_N,\widehat{v}_N)\in\mathcal{B}_N^\delta$, we further linearize \eqref{trun-nonlin-eq} by freezing the advecting velocity at $(u^a+\widehat{u}_N,v^a+\widehat{v}_N)$:
\begin{equation}\label{trun-lin-eq}
\begin{cases}
(u^a+\widehat{u}_N)\partial_xu_N+(v^a+\widehat{v}_N)\partial_yu_N+u_N\partial_xu^a+v_N\partial_yu^a+\partial_xp_N-\epsilon\Delta u_N=-R_u^a,\\
(u^a+\widehat{u}_N)\partial_xv_N+(v^a+\widehat{v}_N)\partial_yv_N+u_N\partial_xv^a+v_N\partial_yv^a+\partial_yp_N-\epsilon\Delta v_N=-R_v^a,\\
\partial_xu_N+\partial_yv_N=0,\\
(u_N,v_N)=(0,0)\text{ on }\partial\Omega_N.
\end{cases}
\end{equation}
It is worth pointing out that this linearization differs from those used in previous work \cite{GI-2023,GN-2017, FGLT-2023,FGLT-2024}, which is essential for establishing the {\it a priori} estimates proving existence.

Assuming that $(u_N,v_N)\in\mathcal{V}_N$ and $p_N\in L^2(\Omega_N)/\mathbb{R}$ solves \eqref{trun-lin-eq} in the weak sense, we now derive uniform {\it a priori} estimates.
\subsection{$L^2$ estimate}\label{L^2}
We start from the $L^2$ estimate, which is one of the most delicate parts of the analysis: the suction dominance in \eqref{alpha0} allows us to get a direct $L^2$ coercivity from the advection term. For later use, we recall the classical Hardy inequality.
\begin{lemma}\label{chardy}
If $f\in H^1(\Omega_N)$ with $f(x,0)=0$, then 
\begin{equation*}
\int_{\Omega_N}\frac{f^2}{y^2}\,{\rm d}x{\rm d}y\leq 4\int_{\Omega_N}|f_y|^2\,{\rm d}x{\rm d}y.
\end{equation*}
\end{lemma}

\begin{proof}
Using $f(x,0)=0$ and the fundamental theorem of calculus, we have for almost every $x$
\begin{equation*}
f(x,y)=\int_0^yf_y(x,s)\,{\rm d}s.
\end{equation*}
It follows from the Cauchy-Schwarz inequality that
\begin{equation*}
f^2(x,y)\leq y\int_0^yf_y^2(x,s)\,{\rm d}s.
\end{equation*}
Then, for any $0<\tau<N$, 
\begin{equation}\label{tau}
\int_{-N}^N\frac{f^2(x,\tau)}{\tau}\,{\rm d}x\leq \int_{-N}^N\int_0^\tau f_y^2(x,s)\,{\rm d}s{\rm d}x\to 0,\quad \text{as $\tau\to 0$}.
\end{equation}

Now integrating by parts in $y$ and applying Young's inequality, we have
\begin{align}
\int_{-N}^N\int_\tau^N\frac{f^2}{y^2}\,{\rm d}y{\rm d}x&=-\int_{-N}^N\frac{f^2(x,N)}{N}\,{\rm d}x+\int_{-N}^N\frac{f^2(x,\tau)}{\tau}\,{\rm d}x+2\int_{-N}^N\int_\tau^N\frac{ff_y}{y}\,{\rm d}y{\rm d}x \nonumber \\
&\leq \int_{-N}^N\frac{f^2(x,\tau)}{\tau}\,{\rm d}x+\frac{1}{2}\int_{-N}^N\int_\tau^N\frac{f^2}{y^2}\,{\rm d}y{\rm d}x+2\int_{-N}^N\int_\tau^N|f_y|^2\,{\rm d}y{\rm d}x. \nonumber
\end{align}
Taking $\tau\to 0$ and using \eqref{tau} in the above inequality completes the proof of Lemma \ref{chardy}.
\end{proof}

Let
\begin{equation*}
\psi_N:=-\int_0^yu_N(x,s)\,{\rm d}s
\end{equation*}
be the stream function of $(u_N,v_N)$. We construct a divergence-free multiplier with the weight $y$:
\begin{equation*}
(m_{N,u},m_{N,v}):=\nabla^\perp(y\psi_N)=(yu_N-\psi_N,yv_N).
\end{equation*}

The following lemma verifies that $(m_{N,u},m_{N,v})$ is an admissible test function for the truncated linearized problem \eqref{trun-lin-eq}.
\begin{lemma}
$(m_{N,u},m_{N,v})\in\mathcal{V}_N$.
\end{lemma}
\begin{proof}
By construction, $(m_{N,u},m_{N,v})\in H^1(\Omega_N)$ is divergence-free, and vanishes on the bottom and two vertical boundaries of $\Omega_N$. At the upper boundary $y=N$,
\begin{equation*}
(\psi_N)_x(x,N)=v_N(x,N)=0,
\end{equation*}
which implies $\psi_N(x,N)$ is constant in $x$. Since $\psi_N(-N,N)=0$, that constant must be zero. This completes the proof.
\end{proof}

Now testing \eqref{trun-lin-eq} with $(m_{N,u},m_{N,v})$, we have
\begin{equation*}
S_1+S_2+S_3=S_4+S_5,
\end{equation*}
where
\begin{align}
S_1&:=\int_{\Omega_N}\Bigl(\bigl((u^a+\widehat{u}_N)\partial_xu_N+(v^a+\widehat{v}_N)\partial_yu_N\bigr)m_{N,u}+\bigl((u^a+\widehat{u}_N)\partial_xv_N+(v^a+\widehat{v}_N)\partial_yv_N\bigr)m_{N,v}\Bigr), \nonumber \\
S_2 &:=\int_{\Omega_N}-\epsilon(\Delta u_Nm_{N,u}+\Delta v_Nm_{N,v}), \nonumber \\
S_3&:=\int_{\Omega_N}(\partial_xp_Nm_{N,u}+\partial_yp_Nm_{N,v}), \nonumber \\
S_4&:=\int_{\Omega_N}-\bigl((u_N\partial_xu^a+v_N\partial_yu^a)m_{N,u}+\bigl(u_N\partial_xv^a+v_N\partial_yv^a\bigr)m_{N,v}\bigr), \nonumber \\
S_5&:=\int_{\Omega_N}-(R_u^am_{N,u}+R_v^am_{N,v}), \nonumber
\end{align}

We estimate term by term as follows.

\paragraph{1. $S_1$}
Substituting in $(m_{N,u},m_{N,v})=y(u_N,v_N)+(-\psi_N,0)$, we split $S_1$ into two parts:
\begin{align}
S_1&=\int_{\Omega_N}y\Bigl(\bigl((u^a+\widehat{u}_N)\partial_xu_N+(v^a+\widehat{v}_N)\partial_yu_N\bigr)u_N+\bigl((u^a+\widehat{u}_N)\partial_xv_N+(v^a+\widehat{v}_N)\partial_yv_N\bigr)v_N\Bigr) \nonumber \\
&-\int_{\Omega_N}\bigl((u^a+\widehat{u}_N)\partial_xu_N+(v^a+\widehat{v}_N)\partial_yu_N\bigr)\psi_N \nonumber \\
&=:\,S_{11}+S_{12}. \nonumber
\end{align}

By integration by parts and incompressibility of $(u^a+\widehat{u}_N,v^a+\widehat{v}_N)$, $S_{11}$ and $S_{12}$ reduce to
\begin{align}
S_{11}&=\int_{\Omega_N}\frac{y}{2}\bigl((u^a+\widehat{u}_N)\partial_x+(v^a+\widehat{v}_N)\partial_y\bigr)(u_N^2+v_N^2) \nonumber \\
&=-\int_{\Omega_N}\frac{y}{2}\bigl(\partial_x(u^a+\widehat{u}_N)+\partial_y(v^a+\widehat{v}_N)\bigr)(u_N^2+v_N^2)-\int_{\Omega_N}\frac{1}{2}(v^a+\widehat{v}_N)(u_N^2+v_N^2) \nonumber \\
&=\int_{\Omega_N}-\frac{1}{2}(v^a+\widehat{v}_N)(u_N^2+v_N^2), \nonumber
\end{align}
and
\begin{align}
S_{12}&=\int_{\Omega_N}u_N\bigl(\partial_x(u^a+\widehat{u}_N)+\partial_y(v^a+\widehat{v}_N)\bigr)\psi_N+\int_{\Omega_N}u_N\bigl((u^a+\widehat{u}_N)\partial_x+(v^a+\widehat{v}_N)\partial_y\bigr)\psi_N \nonumber \\
&=\int_{\Omega_N}u_N\bigl((u^a+\widehat{u}_N)v_N-(v^a+\widehat{v}_N)u_N\bigr). \nonumber
\end{align}

Combining the two parts, $S_1$ reduces to a quadratic form in $(u_N,v_N)$:
\begin{equation*}
S_1=\int_{\Omega_N}\left(-\frac{3}{2}(v^a+\widehat{v}_N)u_N^2+(u^a+\widehat{u}_N)u_Nv_N-\frac{1}{2}(v^a+\widehat{v}_N)v_N^2\right)\geq\theta\|(u_N,v_N)\|_{L^2(\Omega_N)}^2,
\end{equation*}
where
\begin{equation}\label{eval}
\theta:=\operatorname*{ess\,inf}_{\Omega_N}\left(-(v^a+\widehat{v}_N)-\frac{1}{2}\sqrt{(v^a+\widehat{v}_N)^2+(u^a+\widehat{u}_N)^2}\right)
\end{equation}
is the essential infimum over $\Omega_N$ of the smallest eigenvalue of the matrix associated with the quadratic form:
\begin{equation*}
\begin{pmatrix}
-\frac{3}{2}(v^a+\widehat{v}_N)&\frac{1}{2}(u^a+\widehat{u}_N)\\
\frac{1}{2}(u^a+\widehat{u}_N)&-\frac{1}{2}(v^a+\widehat{v}_N)
\end{pmatrix}.
\end{equation*}

\paragraph{2. $S_2$}
We split $S_2$ into two parts:
\begin{equation*}
S_2=-\epsilon\int_{\Omega_N}y\bigl(\Delta u_Nu_N+\Delta v_Nv_N\bigr)+\int_{\Omega_N}\epsilon\psi_N\Delta u_N=:S_{21}+S_{22}.
\end{equation*}

Integration by parts gives
\begin{align}
S_{21}&=\epsilon \int_{\Omega_N}y(|\nabla u_N|^2+|\nabla v_N|^2)+\epsilon\int_{\Omega_N}(\partial_yu_Nu_N+\partial_yv_Nv_N) \nonumber \\
&=\epsilon\int_{\Omega_N}y|\nabla(u_N,v_N)|^2+\frac{\epsilon}{2}\int_{\Omega_N}\partial_y(|u_N|^2+|v_N|^2) \nonumber \\
&=\epsilon\|\sqrt{y}\nabla(u_N,v_N)\|_{L^2(\Omega_N)}^2, \nonumber
\end{align}
and $S_{22}$ vanishes by incompressibility of $(u_N,v_N)$ and its boundary condition:
\begin{align}
S_{22}&=-\epsilon\int_{\Omega_N}\nabla\psi_N\cdot\nabla u_N=-\epsilon\int_{\Omega_N}(v_N\partial_xu_N-u_N\partial_yu_N) \nonumber \\
&=\epsilon\int_{\Omega_N}(v_N\partial_yv_N+u_N\partial_yu_N)=\frac{\epsilon}{2}\int_{\Omega_N}\partial_y(|v_N|^2+|u_N|^2)=0. \nonumber
\end{align}

It follows that
\begin{equation*}
S_2=\epsilon\|\sqrt{y}\nabla(u_N,v_N)\|_{L^2(\Omega_N)}^2.
\end{equation*}

\paragraph{3. $S_3$} $S_3$ vanishes by integration by parts and incompressibility of $(m_{N,u},m_{N,v})$:
\begin{equation*}
S_3=-\int_{\Omega_N}p_N(\partial_xm_{N,u}+\partial_ym_{N,v})=0.
\end{equation*}

\paragraph{4. $S_4$}
We split $S_4$ into two parts:
\begin{align}
S_4&=-\int_{\Omega_N}y\bigl((u_N\partial_xu^a+v_N\partial_yu^a)u_N+(u_N\partial_xv^a+v_N\partial_yv^a)v_N\bigr)+\int_{\Omega_N}\psi_N(u_N\partial_xu^a+v_N\partial_yu^a) \nonumber \\
&=:\,S_{41}+S_{42}. \nonumber
\end{align}

For $S_{41}$, we directly estimate by Hölder's inequality:
\begin{equation*}
|S_{41}|\leq\|y\nabla(u^a,v^a)\|_{L^\infty(\Omega_N)}\|(u_N,v_N)\|_{L^2(\Omega_N)}^2.
\end{equation*}

As for $S_{42}$, it follows from integration by parts and $(u_N,v_N)\perp\nabla\psi_N$ that
\begin{equation*}
S_{42}=-\int_{\Omega_N}u^a(u_N\partial_x\psi_N+v_N\partial_y\psi_N)=0
\end{equation*}
Thus, we have
\begin{equation*}
|S_4|\leq\|y\nabla(u^a,v^a)\|_{L^\infty(\Omega_N)}\|(u_N,v_N)\|_{L^2(\Omega_N)}^2.
\end{equation*}

\paragraph{5. $S_5$} Using H\"{o}lder's inequality and Lemma \ref{chardy}, we estimate $S_5$ as follows:
\begin{align}
S_5&=-\int_{\Omega_N}y(R_u^au_N+R_v^av_N)+\int_{\Omega_N}R_u^a\psi_N \nonumber \\
&\leq\|yR^a\|_{L^2(\Omega_N)}\|(u_N,v_N)\|_{L^2(\Omega_N)}+\|yR^a\|_{L^2(\Omega_N)}\|\psi_N/y\|_{L^2(\Omega_N)} \nonumber \\
&\leq 3\|yR^a\|_{L^2(\Omega_N)}\|(u_N,v_N)\|_{L^2(\Omega_N)} \nonumber \\
&\leq\frac{C_0}{\iota}\|y(R_u^a,R_v^a)\|_{L^2(\Omega_N)}^2+\iota\|(u_N,v_N)\|_{L^2(\Omega_N)}^2, \nonumber
\end{align}

Combining the above estimates, we obtain
\begin{align}
&\bigl(\theta-\|y\nabla(u^a,v^a)\|_{L^\infty(\Omega_N)}-\iota\bigr)\|(u_N,v_N)\|_{L^2(\Omega_N)}^2 \nonumber \\
&\quad+\epsilon\|\sqrt{y}\nabla(u_N,v_N)\|_{L^2(\Omega_N)}^2\leq\frac{C_0}{\iota}\|y(R_u^a,R_v^a)\|_{L^2(\Omega_N)}^2. \label{L^2-est-raw}
\end{align}
From the construction of the approximate profile \eqref{approx-profile} and Proposition \ref{approx-profile-prop}, we have
\begin{equation*}
\begin{gathered}
-(v^a+\widehat{v}_N)=-v_I^0-\widehat{v}_N+O(\epsilon),\quad 
u^a+\widehat{u}_N=u_I^0+u_B^0+\widehat{u}_N+O(\epsilon),\\
y\nabla(u^a,v^a)=y\nabla(u_I^0,v_I^0)+\begin{pmatrix}0&y\partial_yu_B^0\\0&0\end{pmatrix}+O(\epsilon),
\end{gathered}
\end{equation*}
where $O(\epsilon)$ denotes terms bounded by $C\epsilon$ in $L^\infty(\Omega_N)$ with $C>0$ independent of $\epsilon$ and $N$.
Then \eqref{alpha0} and the choice of $(\widehat{u}_N,\widehat{v}_N)$, together with the bounds
\begin{equation*}
|u_B^0|=|\eta e^{-\frac{\alpha y}{\epsilon}}|\leq\|\eta\|_{L^\infty(\mathbb{R})},\qquad|y\partial_yu_B^0|=\left|\eta\frac{\alpha y}{\epsilon}e^{-\frac{\alpha y}{\epsilon}}\right|\leq\|\eta\|_{L^\infty(\mathbb{R})},
\end{equation*}
imply
\begin{align}
&\theta-\|y\nabla(u^a,v^a)\|_{L^\infty(\Omega_N)}-\iota \nonumber \\
&\geq(\lambda_0-\delta)-\frac{1}{2}\sqrt{(\lambda_0-\delta)^2+(\|u_I^0\|_{L^\infty(\Omega_N)}+\|\eta\|_{L^\infty(\mathbb{R})}+\delta)^2} \nonumber \\
&\quad-\|y\nabla(u_I^0,v_I^0)\|_{L^\infty(\Omega_N)}-\|\eta\|_{L^\infty(\mathbb{R})}-\iota-C\epsilon\geq \theta_0, \label{lowebound}
\end{align}
for some constants $C, \theta_0>0$ independent of $\epsilon$ and $N$, whenever $\epsilon$, $\delta$, and $\iota$ are chosen sufficiently small.
Since \eqref{alpha0} gives $\|\eta\|_{L^\infty(\mathbb{R})}<2\lambda_0$, we choose $\theta_0$ so that $\theta_0\leq 1-\|\eta\|_{L^\infty(\mathbb{R})}/(2\lambda_0)$ as well.

Consequently, substituting \eqref{lowebound} into \eqref{L^2-est-raw} and using the residual estimate \eqref{res-est}, we establish the $L^2$ estimate for $(u_N,v_N)$:
\begin{equation}\label{L^2-est}
\|(u_N,v_N)\|_{L^2(\Omega_N)}\leq C\epsilon^2,
\end{equation}
for some constant $C>0$ independent of $\epsilon$ and $N$.

We emphasize that the $L^2$ estimate \eqref{L^2-est} relies on the lower bound \eqref{lowebound}, which is ensured by \eqref{alpha0}.

\subsection{$H^1$ estimate}\label{H^1}
We next derive the $H^1$ estimate, another particularly delicate step in the analysis. Exploiting the explicit expression for the leading boundary-layer profile $u_B^0$, we establish the following Hardy-type inequalities, which provide a small coefficient needed to absorb the singular term involving $\partial_yu_B^0$ into the diffusion term.

\begin{lemma}\label{hardy}
Define the boundary-layer weight
\begin{equation*}
w=|\partial_y u_B^0|=\frac{|\eta|\alpha}{\epsilon}e^{-\alpha y/\epsilon}.
\end{equation*}
If $f\in H^1(\Omega)$ with $f(x,0)=0$, then
\begin{align}
\int_\Omega w|f|^2\,{\rm d}x{\rm d}y&\leq\epsilon\frac{\|\eta\|_{L^\infty(\mathbb{R})}}{\lambda_0}\|f_y\|_{L^2(\Omega)}^2, \label{hardy1} \\
\int_\Omega w^2|f|^2\,{\rm d}x{\rm d}y&\leq\frac{\|\eta\|_{L^\infty(\mathbb{R})}^2}{4}\|f_y\|_{L^2(\Omega)}^2. \label{hardy2}
\end{align}

The same inequalities hold with $\Omega$ replaced by $\Omega_N$
for $f\in H_0^1(\Omega_N)$, by zero-extending $f$ to $\Omega$.
\end{lemma}
\begin{proof}
The same argument as in the proof of Lemma \ref{chardy} gives
\begin{equation*}
|f(x,y)|^2\leq y\int_0^y|f_y(x,s)|^2\,{\rm d}s.
\end{equation*}
Since $\alpha\geq\lambda_0$, we have
\begin{align}
\int_0^\infty w(x,y)|f(x,y)|^2\,{\rm d}y&\leq\int_0^\infty|f_y(x,s)|^2\left(\int_0^\infty|\eta(x)|\frac{\alpha(x)}{\epsilon}ye^{-\alpha(x)y/\epsilon}\,{\rm d}y\right)\,{\rm d}s \nonumber \\
&\leq\|\eta\|_{L^\infty(\mathbb{R})}\frac{\epsilon}{\alpha(x)}\int_0^\infty|f_y(x,s)|^2\,{\rm d}s\leq\epsilon\frac{\|\eta\|_{L^\infty(\mathbb{R})}}{\lambda_0}\int_0^\infty|f_y(x,s)|^2\,{\rm d}s. \nonumber
\end{align}
Integrating in $x$ proves \eqref{hardy1}.

Similarly, we have
\begin{align}
\int_0^\infty w^2(x,y)|f(x,y)|^2\,{\rm d}y&\leq\int_0^\infty|f_y(x,s)|^2\left(\int_s^\infty|\eta(x)|^2\frac{\alpha^2(x)}{\epsilon^2}ye^{-2\alpha(x)y/\epsilon}\,{\rm d}y\right)\,{\rm d}s \nonumber \\
&\leq\frac{\|\eta\|_{L^\infty(\mathbb{R})}^2}{4}\int_0^\infty|f_y(x,s)|^2\,{\rm d}s. \nonumber
\end{align}
Integration in $x$ proves \eqref{hardy2}.
\end{proof}

We now
test \eqref{trun-lin-eq} with $(u_N,v_N)$ to see
\begin{equation*}
Q_1+Q_2+Q_3=Q_4+Q_5,
\end{equation*}
where
\begin{align}
Q_1&:=\int_{\Omega_N}-\epsilon(\Delta u_Nu_N+\Delta v_Nv_N), \nonumber \\
Q_2&:=\int_{\Omega_N}\Bigl(\bigl((u^a+\widehat{u}_N)\partial_xu_N+(v^a+\widehat{v}_N)\partial_yu_N\bigr)u_N+\bigl((u^a+\widehat{u}_N)\partial_xv_N+(v^a+\widehat{v}_N)\partial_yv_N\bigr)v_N\Bigr), \nonumber \\
Q_3&:=\int_{\Omega_N}(\partial_xp_Nu_N+\partial_yp_Nv_N), \nonumber \\
Q_4&:=\int_{\Omega_N}-\bigl((u_N\partial_xu^a+v_N\partial_yu^a)u_N+\bigl(u_N\partial_xv^a+v_N\partial_yv^a\bigr)v_N\bigr), \nonumber \\
Q_5&:=\int_{\Omega_N}-(R_u^au_N+R_v^av_N). \nonumber
\end{align}

We estimate term by term as follows.

\paragraph{1. $Q_1$}
A direct integration by parts gives
\begin{equation*}
Q_1=\epsilon\int_{\Omega_N}|\nabla(u_N,v_N)|^2.
\end{equation*}

\paragraph{2. $Q_2$}
$Q_2$ vanishes by integration by parts and incompressibility of $(u^a+\widehat{u}_N,v^a+\widehat{v}_N)$:
\begin{align}
Q_2&=\int_{\Omega_N}\frac{1}{2}\bigl((u^a+\widehat{u}_N)\partial_x+(v^a+\widehat{v}_N)\partial_y\bigr)(u_N^2+v_N^2) \nonumber \\
&=-\frac{1}{2}\int_{\Omega_N}\bigl(\partial_x(u^a+\widehat{u}_N)+\partial_y(v^a+\widehat{v}_N)\bigr)(u_N^2+v_N^2)=0. \nonumber
\end{align}

\paragraph{3. $Q_3$}
Similarly, $Q_3$ vanishes by integration by parts and incompressibility of $(u_N,v_N)$:
\begin{equation*}
Q_3=-\int_{\Omega_N}p_N(\partial_xu_N+\partial_yv_N)=0.
\end{equation*}

\paragraph{4. $Q_4$}
We separate the most singular term $\partial_yu_B^0$ from the rest of the terms in $Q_4$:
\begin{align}
Q_4&=-\int_{\Omega_N}v_N\partial_yu_B^0u_N-\int_{\Omega_N}\bigl(u_N\partial_xu^au_N+v_N\partial_y(u^a-u_B^0)u_N+u_N\partial_xv^av_N+v_N\partial_yv^av_N\bigr) \nonumber \\
&=:\,Q_{41}+Q_{42}. \nonumber
\end{align}

For $Q_{41}$, we invoke the Hardy-type inequality \eqref{hardy1} to obtain
\begin{equation*}
|Q_{41}|\leq\int_{\Omega_N}w|u_N||v_N|\leq\frac{1}{2}\int_{\Omega_N}w(|u_N|^2+|v_N|^2)\leq\frac{\epsilon\|\eta\|_{L^\infty(\mathbb{R})}}{2\lambda_0}\|\partial_y(u_N,v_N)\|_{L^2(\Omega_N)}^2.
\end{equation*}

By construction of the approximate solution \eqref{approx-profile} and Proposition \ref{approx-profile-prop}, each term in the expansions of $\partial_xu^a$, $\partial_y(u^a-u_B^0)$, $\partial_xv^a$, and $\partial_yv^a$ is bounded in $L^\infty(\Omega_N)$ uniformly in $\epsilon$ and $N$. Therefore, $Q_{42}$ can be controlled directly by Hölder's inequality:
\begin{align}
|Q_{42}|&\leq\|\partial_xu^a\|_{L^\infty(\Omega_N)}\|u_N\|_{L^2(\Omega_N)}^2+\|\partial_y(u^a-u_B^0)\|_{L^\infty(\Omega_N)}\|v_N\|_{L^2(\Omega_N)}\|u_N\|_{L^2(\Omega_N)} \nonumber \\
&\quad+\|\partial_xv^a\|_{L^\infty(\Omega_N)}\|u_N\|_{L^2(\Omega_N)}\|v_N\|_{L^2(\Omega_N)}+\|\partial_yv^a\|_{L^\infty(\Omega_N)}\|v_N\|_{L^2(\Omega_N)}^2 \nonumber \\
&\leq C\|(u_N,v_N)\|_{L^2(\Omega_N)}^2, \nonumber
\end{align}
for some constant $C>0$ independent of $\epsilon$ and $N$.

Combining the two parts yields
\begin{equation*}
|Q_4|\leq\frac{\epsilon\|\eta\|_{L^\infty(\mathbb{R})}}{2\lambda_0}\|\partial_y(u_N,v_N)\|_{L^2(\Omega_N)}^2+C\|(u_N,v_N)\|_{L^2(\Omega_N)}^2.
\end{equation*}

\paragraph{5. $Q_5$} A direct use of H\"{o}lder's inequality gives
\begin{equation*}
|Q_5|\leq\|(R_u^a,R_v^a)\|_{L^2(\Omega_N)}\|(u_N,v_N)\|_{L^2(\Omega_N)}.
\end{equation*}

Collecting the above estimates, we have
\begin{align}
&\epsilon\Big(1-\frac{\|\eta\|_{L^\infty(\mathbb{R})}}{2\lambda_0}\Big)\|\nabla(u_N,v_N)\|_{L^2(\Omega_N)}^2 \nonumber \\
&\leq C\|(u_N,v_N)\|_{L^2(\Omega_N)}^2+\|(R_u^a,R_v^a)\|_{L^2(\Omega_N)}\|(u_N,v_N)\|_{L^2(\Omega_N)}. \label{H^1-est-raw}
\end{align}

The inequality \eqref{alpha0} in particular implies $\|\eta\|_{L^\infty(\mathbb{R})}<2\lambda_0$. Together with the residual estimate \eqref{res-est} and the previously established $L^2$ estimate \eqref{L^2-est}, we conclude the $H^1$ estimate for $(u_N,v_N)$:
\begin{equation}\label{H^1-est}
\|(u_N,v_N)\|_{H^1(\Omega_N)}\leq C\epsilon^{5/4},
\end{equation}
for some constant $C>0$ independent of $\epsilon$ and $N$.

\subsection{$H^2$ estimate}\label{H^2}
Lastly, we establish the $H^2$ estimates via the Stokes estimates. To ensure the estimate is uniform in $N$, attention must be paid to the dependence of the Stokes estimate on the domain.

We first introduce the classical Stokes estimate on convex polygons.
\begin{lemma}[Kellogg--Osborn {\cite[Theorem~2]{KO-1976}}]\label{stokes-polygon}
Let $P\subset\mathbb{R}^2$ be a bounded convex polygon and $(F_u,F_v)\in L^2(P)$. If $(u,v)\in H^1_0(P)$ and $p\in L^2(P)/\mathbb{R}$ solve the Stokes system
\begin{equation}\label{stokes-system}
\begin{cases}
-\Delta(u,v)+\nabla p=(F_u,F_v),\\
\partial_xu+\partial_yv=0,
\end{cases}
\end{equation}
then we have $(u,v)\in H^2(P)$ and $p\in H^1(P)/\mathbb{R}$, with the estimate
\begin{equation*}
\|(u,v)\|_{H^2(P)}+\|\nabla p\|_{L^2(P)}\leq C\|(F_u,F_v)\|_{L^2(P)},
\end{equation*}
for some constant $C>0$ depending only on $P$. Here $L^2(P)/\mathbb{R}$ and $H^1(P)/\mathbb{R}$ denote the corresponding spaces modulo constants.
\end{lemma}

The corresponding dilation-invariant version is as follows.
\begin{corollary}[Uniform Stokes Estimate]\label{stokes-polygon-uniform}
Let $P\subset\mathbb{R}^2$ be a bounded convex polygon and $P_\lambda=\lambda P$ for $\lambda>0$. If $(u_\lambda,v_\lambda)\in H_0^1(P_\lambda)$ and $p_\lambda\in L^2(P_\lambda)/\mathbb{R}$ solve \eqref{stokes-system} on $P_\lambda$ with forcing $(F_{\lambda,u},F_{\lambda,v})\in L^2(P_\lambda)$, then we have
\begin{equation*}
\|\nabla^2(u_\lambda,v_\lambda)\|_{L^2(P_\lambda)}+\|\nabla p_\lambda\|_{L^2(P_\lambda)}\leq C\|(F_{\lambda,u},F_{\lambda,v})\|_{L^2(P_\lambda)},
\end{equation*}
for some constant $C>0$ independent of $\lambda$.
\end{corollary}
\begin{proof}
Set $(x,y)=\lambda(\widetilde{x},\widetilde{y})$ and define on $P$
\begin{align}
(\widetilde{u},\widetilde{v})(\widetilde{x},\widetilde{y})&=(u_\lambda,v_\lambda)(\lambda\widetilde{x},\lambda\widetilde{y}), \nonumber \\
\widetilde{p}(\widetilde{x},\widetilde{y})&=\lambda p_\lambda(\lambda\widetilde{x},\lambda\widetilde{y}), \nonumber \\
(\widetilde{F}_u,\widetilde{F}_v)(\widetilde{x},\widetilde{y})&=\lambda^2(F_{\lambda,u},F_{\lambda,v})(\lambda\widetilde{x},\lambda\widetilde{y}). \nonumber
\end{align}
Then $(\widetilde{u},\widetilde{v})\in H^1_0(P)$ and $\widetilde{p}\in L^2(P)/\mathbb{R}$ solve
\begin{equation*}
\begin{cases}
-\Delta_{\widetilde{x},\widetilde{y}}(\widetilde{u},\widetilde{v})+\nabla_{\widetilde{x},\widetilde{y}}\widetilde{p}=(\widetilde{F}_u,\widetilde{F}_v),\\
\partial_{\widetilde{x}}\widetilde{u}+\partial_{\widetilde{y}}\widetilde{v}=0.
\end{cases}
\end{equation*}
By Lemma \ref{stokes-polygon},
\begin{equation*}
\|\nabla_{\widetilde{x},\widetilde{y}}^2(\widetilde{u},\widetilde{v})\|_{L^2(P)}+\|\nabla_{\widetilde{x},\widetilde{y}}\widetilde{p}\|_{L^2(P)}\leq C\|(\widetilde{F}_u,\widetilde{F}_v)\|_{L^2(P)}.
\end{equation*}
In two dimensions, we have the following scaling relations:
\begin{align}
\|\nabla_{\widetilde{x},\widetilde{y}}^2(\widetilde{u},\widetilde{v})\|_{L^2(P)}&=\lambda\|\nabla_{x,y}^2(u_\lambda,v_\lambda)\|_{L^2(P_\lambda)}, \nonumber \\
\|\nabla_{\widetilde{x},\widetilde{y}}\widetilde{p}\|_{L^2(P)}&=\lambda\|\nabla_{x,y}p_\lambda\|_{L^2(P_\lambda)}, \nonumber \\
\|(\widetilde{F}_u,\widetilde{F}_v)\|_{L^2(P)}&=\lambda\|(F_{\lambda,u},F_{\lambda,v})\|_{L^2(P_\lambda)}. \nonumber
\end{align}
Substituting these relations and canceling $\lambda$ proves the desired estimate.
\end{proof}

For the uniform $L^\infty$ estimate, we also provide the following dilation-invariant Agmon inequality.
\begin{corollary}\label{agmon-uniform}
Let $D\subset\mathbb{R}^2$ be a bounded Sobolev extension domain and $D_\lambda=\lambda D$ be a dilation of $D$ by a factor $\lambda>0$. Then, for $f_\lambda\in H^2(D_\lambda)\cap H_0^1(D_\lambda)$, we have
\begin{equation*}
\|f_\lambda\|_{L^\infty(D_\lambda)}\leq C\|f_\lambda\|_{L^2(D_\lambda)}^{1/2}\|\nabla^2f_\lambda\|_{L^2(D_\lambda)}^{1/2},
\end{equation*}
for some constant $C>0$ independent of $\lambda$.
\end{corollary}
\begin{proof}
For a fixed domain $D$, the usual Agmon inequality reads
\begin{equation}\label{Agmon}
\|f\|_{L^\infty(D)}\leq C\|f\|_{L^2(D)}^{1/2}\|\nabla^2f\|_{L^2(D)}^{1/2}.
\end{equation}
Set $\xi=\lambda\widetilde{\xi}$ and define $\widetilde{f}(\widetilde{\xi})=f_\lambda(\lambda\widetilde{\xi})$ on $D$. In two dimensions, the following scaling relations hold:
\begin{equation*}
\|\widetilde{f}\|_{L^\infty(D)}=\|f_\lambda\|_{L^\infty(D_\lambda)},\quad\|\widetilde{f}\|_{L^2(D)}=\lambda^{-1}\|f_\lambda\|_{L^2(D_\lambda)},\quad\|\nabla_{\widetilde{\xi}}^2\widetilde{f}\|_{L^2(D)}=\lambda\|\nabla_\xi^2f_\lambda\|_{L^2(D_\lambda)}.
\end{equation*}
Substituting these relations into \eqref{Agmon} completes the proof of Corollary \ref{agmon-uniform}.
\end{proof}

Now we are in a position to establish the $H^2$ estimate.

Rewrite \eqref{trun-lin-eq} as
\begin{equation*}
\begin{cases}
-\epsilon\Delta u_N+\partial_xp_N=F_u,\\
-\epsilon\Delta v_N+\partial_yp_N=F_v,
\end{cases}
\end{equation*}
where the source terms $F_u$ and $F_v$ are given by
\begin{align}
F_u&:=-\bigl((u^a+\widehat{u}_N)\partial_xu_N+(v^a+\widehat{v}_N)\partial_yu_N+u_N\partial_xu^a+v_N\partial_y(u^a-u_B^0)\bigr)-v_N\partial_yu_B^0-R_u^a \nonumber \\
&=:F_{u,1}-v_N\partial_yu_B^0-R_u^a, \nonumber \\
F_v&:=-\bigl((u^a+\widehat{u}_N)\partial_xv_N+(v^a+\widehat{v}_N)\partial_yv_N+u_N\partial_xv^a+v_N\partial_yv^a\bigr)-R_v^a \nonumber \\
&=:F_{v,1}-R_v^a. \nonumber
\end{align}

As mentioned in the estimate of $Q_4$ in Section \ref{H^1}, $u^a+\widehat{u}_N$, $v^a+\widehat{v}_N$, $\partial_xu^a$, $\partial_xv^a$, $\partial_y(u^a-u_B^0)$, and $\partial_yv^a$ are bounded in $L^\infty(\Omega_N)$ uniformly in $\epsilon$ and $N$, hence $F_{u,1}$ and $F_{v,1}$ are directly controlled as follows:
\begin{equation*}
\|(F_{u,1},F_{v,1})\|_{L^2(\Omega_N)}\leq C\|(u_N,v_N)\|_{H^1(\Omega_N)},
\end{equation*}
for some constant $C>0$ independent of $\epsilon$ and $N$.

The singular term $v_N\partial_yu_B^0$ is handled by the Hardy-type inequality \eqref{hardy2}:
\begin{equation*}
\|v_N\partial_yu_B^0\|_{L^2(\Omega_N)}=\|wv_N\|_{L^2(\Omega_N)}\leq\frac{\|\eta\|_{L^\infty(\mathbb{R})}}{2}\|\partial_yv_N\|_{L^2(\Omega_N)}.
\end{equation*}

Combining the above estimates, then applying the residual estimate \eqref{res-est} and the previously established $H^1$ estimate \eqref{H^1-est}, we obtain
\begin{equation*}
\|(F_u,F_v)\|_{L^2(\Omega_N)}\leq C\|(u_N,v_N)\|_{H^1(\Omega_N)}+\|(R_u^a,R_v^a)\|_{L^2(\Omega_N)}\leq C\epsilon^{5/4},
\end{equation*}
for some constant $C>0$ independent of $\epsilon$ and $N$. Consequently, the uniform Stokes estimate, Corollary \ref{stokes-polygon-uniform}, yields
\begin{equation*}
\epsilon\|\nabla^2(u_N,v_N)\|_{L^2(\Omega_N)}+\|\nabla p_N\|_{L^2(\Omega_N)}\leq C\|(F_u,F_v)\|_{L^2(\Omega_N)}\leq C\epsilon^{5/4},
\end{equation*}
and thus we conclude
\begin{equation}\label{H^2-est}
\|\nabla^2(u_N,v_N)\|_{L^2(\Omega_N)}\leq C\epsilon^{1/4},\qquad\|\nabla p_N\|_{L^2(\Omega_N)}\leq C\epsilon^{5/4},
\end{equation}
for some constant $C>0$ independent of $\epsilon$ and $N$.

{Finally, since $(u_N,v_N)\in H^2(\Omega_N)\cap H^1_0(\Omega_N)$,} the $L^2$ estimate \eqref{L^2-est} and the uniform Agmon inequality, Corollary \ref{agmon-uniform}, give the $L^\infty$ estimate
\begin{equation}\label{L^infty-est}
\|(u_N,v_N)\|_{L^\infty(\Omega_N)}\leq C\|(u_N,v_N)\|_{L^2(\Omega_N)}^{1/2}\|\nabla^2(u_N,v_N)\|_{L^2(\Omega_N)}^{1/2}\leq C\epsilon^{9/8},
\end{equation}
for some constant $C>0$ independent of $\epsilon$ and $N$. In fact, since $H^2(\Omega_N)\hookrightarrow C(\overline{\Omega_N})$, we have
\begin{equation*}
\|(u_N,v_N)\|_{C(\overline{\Omega_N})}=\|(u_N,v_N)\|_{L^\infty(\Omega_N)}\leq C\epsilon^{9/8},
\end{equation*}
for some constant $C>0$ independent of $\epsilon$ and $N$.

\section{Existence of the Remainder System}\label{exist}
In this section, we subsequently establish the existence of a solution to the truncated linear problem \eqref{trun-lin-eq} and the original global nonlinear problem \eqref{remainder-eq}.

For simplicity, denote $U_N=(u_N,v_N)$, $\widehat{U}_N=(\widehat{u}_N,\widehat{v}_N)$, $U=(u,v)$, $U^a=(u^a,v^a)$, and $R^a=(R_u^a,R_v^a)$. Then \eqref{remainder-eq}, \eqref{trun-nonlin-eq}, and \eqref{trun-lin-eq} read, respectively.
\begin{align}
&\begin{cases}
\bigl((U^a+U)\cdot\nabla\bigr)U+(U\cdot\nabla)U^a+\nabla p-\epsilon\Delta U=-R^a,\\
\nabla\cdot U=0,\\
U(x,0)=0,\\
\lim\limits_{y\to\infty}U(x,y)=0,
\end{cases} \label{remainder-eq-vec} \\
&\begin{cases}
\bigl((U^a+U_N)\cdot\nabla\bigr)U_N+(U_N\cdot\nabla)U^a+\nabla p_N-\epsilon\Delta U_N=-R^a,\\
\nabla\cdot U_N=0,\\
U_N=0\text{ on }\partial\Omega_N,
\end{cases} \label{trun-nonlin-eq-vec}
\end{align}

and
\begin{equation}\label{trun-lin-eq-vec}
\begin{cases}
\bigl((U^a+\widehat{U}_N)\cdot\nabla\bigr)U_N+(U_N\cdot\nabla)U^a+\nabla p_N-\epsilon\Delta U_N=-R^a,\\
\nabla\cdot U_N=0,\\
U_N=0\text{ on }\partial\Omega_N.
\end{cases}
\end{equation}

\subsection{Existence of truncated linear problem}\label{trun-lin-eq-sec}
We first prove the existence and uniqueness for the truncated linear problem \eqref{trun-lin-eq-vec} in $\mathcal{V}_N$ by the Lax--Milgram theorem followed by the Fredholm alternative, thus justifying the {\it a priori} estimates established in Section \ref{est}.

\begin{proposition}\label{trun-lin-eq-exist}
For sufficiently small $\delta,\epsilon>0$, \eqref{trun-lin-eq-vec} admits a unique strong solution $(U_N,p_N)\in(\mathcal{V}_N\cap H^2(\Omega_N))\times(H^1(\Omega_N)/\mathbb{R})$ satisfying the estimates:
\begin{equation}\label{H2est}
\begin{gathered}
\|U_N\|_{L^2(\Omega_N)}\leq C\epsilon^2,\qquad\|\nabla U_N\|_{L^2(\Omega_N)}\leq C\epsilon^{5/4},\qquad\|\nabla^2U_N\|_{L^2(\Omega_N)}\leq C\epsilon^{1/4},\\
\|U_N\|_{C(\overline{\Omega_N})}\leq C\epsilon^{9/8},\qquad\|\nabla p_N\|_{L^2(\Omega_N)}\leq C\epsilon^{5/4},
\end{gathered}
\end{equation}
for some constant $C>0$ independent of $\epsilon$ and $N$.
\end{proposition}

\begin{proof}
We proceed in four steps.

\noindent\textit{1}.
For $U,\Phi\in\mathcal{V}_N$, we define the bilinear form
\begin{equation*}
a_{\widehat{U}_N}(U,\Phi):=\epsilon\int_{\Omega_N}\nabla U:\nabla\Phi+\int_{\Omega_N}\bigl((U^a+\widehat{U}_N)\cdot\nabla U\bigr)\cdot\Phi+\int_{\Omega_N}(U\cdot\nabla U^a)\cdot\Phi,
\end{equation*}
and the linear functional
\begin{equation*}
\ell(\Phi):=-\int_{\Omega_N}R^a\cdot\Phi.
\end{equation*}

The weak formulation of \eqref{trun-lin-eq-vec} is equivalent to finding $U\in\mathcal{V}_N$ such that
\begin{equation}\label{trun-lin-eq-wk}
a_{\widehat{U}_N}(U,\Phi)=\ell(\Phi),\quad\forall\Phi\in\mathcal{V}_N.
\end{equation}

\noindent\textit{2}.
The argument leading to the $H^1$ estimate \eqref{H^1-est-raw} gives the G{\aa}rding inequality:
\begin{equation*}
a_{\widehat{U}_N}(U,U)\geq\epsilon\theta_0\|\nabla U\|_2^2-C\|U\|_2^2,\quad\forall U\in\mathcal{V}_N.
\end{equation*}
To obtain coercivity, we choose a constant $\Lambda>C$ and define the bilinear form
\begin{equation*}
a_{\widehat{U}_N}^\Lambda(U,\Phi):=a_{\widehat{U}_N}(U,\Phi)+\Lambda\int_{\Omega_N}U\cdot\Phi.
\end{equation*}
Then, $a_{\widehat{U}_N}^\Lambda$ is coercive:
\begin{equation*}
a_{\widehat{U}_N}^\Lambda(U,U)\geq\epsilon\theta_0\|\nabla U\|_2^2+(\Lambda-C)\|U\|_2^2\geq c_{\widehat{U}_N}^\Lambda\|U\|_{H^1}^2,
\end{equation*}
for some constant $c_{\widehat{U}_N}^\Lambda>0$. Moreover, $a_{\widehat{U}_N}^\Lambda$ and $\ell$ are continuous:
\begin{equation*}
\begin{gathered}
\big|a_{\widehat{U}_N}^\Lambda(U,\Phi)\big|\leq(\epsilon+\|U^a\|_\infty+\|\widehat{U}_N\|_\infty+\|\nabla U^a\|_\infty+\Lambda)\|U\|_{H^1}\|\Phi\|_{H^1},\\
|\ell(\Phi)|\leq\|R^a\|_2\|\Phi\|_2\leq\|R^a\|_2\|\Phi\|_{H^1}.
\end{gathered}
\end{equation*}
The Lax--Milgram theorem then provides a unique solution $U\in\mathcal{V}_N$ to the weak problem:
\begin{equation*}
a_{\widehat{U}_N}^{\Lambda}(U,\Phi)=\ell(\Phi),\quad\forall\Phi\in\mathcal{V}_N.
\end{equation*}

\noindent\textit{3}.
For $U,\Phi\in\mathcal{V}_N$, define the operators $A_{\widehat{U}_N},J,A_{\widehat{U}_N}^\Lambda:\mathcal{V}_N\to\mathcal{V}_N^*$ by
\begin{equation*}
\big\langle A_{\widehat{U}_N}U,\Phi\big\rangle:=a_{\widehat{U}_N}(U,\Phi),\quad 
\langle JU,\Phi\rangle:=\int_{\Omega_N}U\cdot\Phi,\quad 
A_{\widehat{U}_N}^\Lambda:=A_{\widehat{U}_N}+\Lambda J.
\end{equation*}
Then the weak formulation \eqref{trun-lin-eq-wk} can be written as $A_{\widehat{U}_N}U=\ell$, or equivalently
\begin{equation*}
\big(A_{\widehat{U}_N}^\Lambda-\Lambda J\big)U=\ell.
\end{equation*}
By step 2, $A_{\widehat{U}_N}^\Lambda$ is invertible, so the equation can be further written as
\begin{equation*}
\big(I-\Lambda\big(A_{\widehat{U}_N}^\Lambda\big)^{-1}J\big)U=\big(A_{\widehat{U}_N}^\Lambda\big)^{-1}\ell.
\end{equation*}

Since $\mathcal{V}_N\hookrightarrow L^2(\Omega_N)$ is compact and $L^2(\Omega_N)\rightarrow\mathcal{V}_N^*$ is continuous, $J:\mathcal{V}_N\hookrightarrow\mathcal{V}_N^*$ is compact. Moreover, the coercivity of $a_{\widehat{U}_N}^\Lambda$ implies $\big(A_{\widehat{U}_N}^\Lambda\big)^{-1}$ is continuous, so that $K_{\widehat{U}_N}^\Lambda:=\Lambda\big(A_{\widehat{U}_N}^\Lambda\big)^{-1}J$ is compact. On the other hand, it follows from the $L^2$ estimates \eqref{L^2-est-raw} with $R^a=0$ that 
$$\operatorname{Ker}\big(I-K_{\widehat{U}_N}^\Lambda\big)=\{0\}.$$
Thus, it follows from the Fredholm alternative that \eqref{trun-lin-eq-wk} admits a unique solution $U_N\in \mathcal{V}_{N}$.

\noindent\textit{4}.
Define
\begin{equation*}
\mathcal{F}_N(\Phi):=a_{\widehat{U}_N}(U_N,\Phi)-\ell(\Phi),\quad\text{for }\Phi\in H_0^1(\Omega_N).
\end{equation*}
By step 2, $\mathcal{F}_N$ is continuous on $H_0^1(\Omega_N)$ and vanishes on $\mathcal{V}_N$ since $U_N$ solves \eqref{trun-lin-eq-wk}. Then, it follows from \cite[Corollary~III.5.1]{Galdi-2011} that there exists a pressure $p_N\in L^2(\Omega_N)$, unique up to an additive constant, such that
\begin{equation*}
\mathcal{F}_N=-\nabla p_N.
\end{equation*}

Consequently, $(U_N,p_N)\in\mathcal{V}_N\times(L^2(\Omega_N)/\mathbb{R})$ solves \eqref{trun-lin-eq-vec}, and the {\it a priori} estimates established in Section \ref{est} further imply $(U_N,p_N)\in(\mathcal{V}_{N}\cap H^2(\Omega_{N}))\times(H^1(\Omega_{N})/\mathbb{R})$ and satisfies the estimates \eqref{H2est}. This completes the proof of Proposition \ref{trun-lin-eq-exist}.
\end{proof}

\subsection{Existence of global nonlinear problem}\label{remainder-eq-sec}
We now prove the existence and uniform estimates for the global nonlinear problem \eqref{remainder-eq-vec}. We first apply Schauder's fixed point theorem to obtain a solution to the truncated nonlinear problem \eqref{trun-nonlin-eq-vec} and then pass to the limit $N\to \infty$ to obtain a solution to the original global nonlinear problem \eqref{remainder-eq-vec}.

\begin{proposition}\label{global-nonlin-eq-exist}
For sufficiently small $\delta,\epsilon>0$, with $\epsilon$ sufficiently small relative to $\delta$, \eqref{remainder-eq-vec} admits a strong solution $(U,p)\in(H_0^1(\Omega)\cap H^2(\Omega))\times(H_{\mathrm{loc}}^1(\Omega)/\mathbb{R})$ satisfying the estimates:
\begin{equation*}
\begin{gathered}
\|U\|_{L^2(\Omega)}\leq C\epsilon^2,\qquad\|\nabla U\|_{L^2(\Omega)}\leq C\epsilon^{5/4},\qquad\|\nabla^2U\|_{L^2(\Omega)}\leq C\epsilon^{1/4},\\
\|U\|_{C(\overline{\Omega})}\leq C\epsilon^{9/8},\qquad\|\nabla p\|_{L^2(\Omega)}\leq C\epsilon^{5/4},
\end{gathered}
\end{equation*}
for some constant $C>0$ independent of $\epsilon$.
\end{proposition}

\begin{proof}
We proceed in six steps.

\noindent\textit{1}. For sufficiently small $\delta,\epsilon>0$, let $T_N:\mathcal{B}_N^\delta\to C(\overline{\Omega_N})$ be the solution map $\widehat{U}_N\mapsto U_N$ of the truncated linear problem \eqref{trun-lin-eq-vec}. By Proposition \ref{trun-lin-eq-exist}, $T_N$ is well-defined and satisfies the uniform-in-$N$ estimate:
\begin{equation*}
\|T_N(\widehat{U}_N)\|_{L^\infty(\Omega_N)}\leq C\epsilon^{9/8},
\end{equation*}
which ensures $T_N$ maps $\mathcal{B}_N^\delta$ into itself whenever $C\epsilon^{9/8}\leq\delta$. Moreover, since $T_N(\mathcal{B}_N^\delta)$ is bounded in $H^2(\Omega_N)$ by Proposition \ref{trun-lin-eq-exist} and $H^2(\Omega_N)\hookrightarrow C(\overline{\Omega_N})$ is compact, $T_N$ is compact.

On the other hand, for any sequence $\widehat{U}_N^j\to\widehat{U}_N^\ast$ in $\mathcal{B}_N^\delta$, the boundedness of $U_N^j:=T_N(\widehat{U}_N^j)$ in $H^2(\Omega_N)$ together with the compact embedding $H^2(\Omega_N)\hookrightarrow H^1(\Omega_N)$ and $C(\overline{\Omega_N})$ yields, up to a subsequence,
\begin{equation*}
U_N^j\to U_N^\ast,\quad\text{in }H^1(\Omega_N)\cap C(\overline{\Omega_N}),
\end{equation*}
for some $U_N^\ast\in H^1(\Omega_N)\cap C(\overline{\Omega_N})$. Taking the limit $j\to \infty$ in the weak formulation \eqref{trun-lin-eq-wk}, we have
\begin{equation*}
a_{\widehat{U}_N^\ast}(U_N^\ast,\Phi)=\ell(\Phi),\quad\forall\Phi\in\mathcal{V}_N.
\end{equation*}
The uniqueness then indicates $U_N^\ast=T_N(\widehat{U}_N^\ast)$. The same argument shows that every subsequence of $\{U_N^j\}$ has a further subsequence converging to $T_N(\widehat{U}_N^\ast)$. Hence the full sequence converges to this limit, proving the continuity of $T_N$.

Therefore, by Schauder's fixed point theorem \cite[Corollary~11.2]{GT-2001}, there exists a $U_{N}=T_N(U_N)\in\mathcal{V}_N\cap H^2(\Omega_N)$, which, together with its associated pressure $p_N\in H^1(\Omega_{N})/\mathbb{R}$, solves the truncated nonlinear problem \eqref{trun-nonlin-eq-vec} and satisfies the estimates in Proposition \ref{trun-lin-eq-exist}.

\noindent\textit{2}. For each $N$, extend $U_N\in H^1_0(\Omega_N)$ to $H^1_0(\Omega)$ by zero and denote the extension by the same symbol. Since $U_N$ is bounded by $C\epsilon^{5/4}$ in $H^1_0(\Omega)$, up to a subsequence, we have
\begin{equation*}
U_N\rightharpoonup U,\quad\text{in }H^1(\Omega),
\end{equation*}
for some $U\in H_0^1(\Omega)$, and lower semicontinuity gives
\begin{equation*}
\|U\|_{L^2(\Omega)}\leq C\epsilon^2,\qquad\|\nabla U\|_{L^2(\Omega)}\leq C\epsilon^{5/4}.
\end{equation*}

\noindent\textit{3}. Since $\|U_N\|_{H^2(\Omega_{N})}\leq C\epsilon^{1/4}$ is uniformly bounded in $N$, a diagonal argument yields, along a further subsequence,
\begin{equation*}
U_N\rightharpoonup U,\quad\text{in }H^2_{\mathrm{loc}}(\Omega),
\end{equation*}
and lower semicontinuity gives, for each $N$,
\begin{equation*}
\|\nabla^2U\|_{L^2(\Omega_N)}\leq C\epsilon^{1/4}.
\end{equation*}
Sending $N$ to infinity yields that $U\in H^2(\Omega)$ with
\begin{equation*}
\|\nabla^2U\|_{L^2(\Omega)}\leq C\epsilon^{1/4}.
\end{equation*}

\noindent\textit{4}.
The compact embedding $H^2(\Omega_N)\hookrightarrow H^1(\Omega_N)$ and $C(\overline{\Omega_N})$ for each $N$ together with a diagonal argument yields, along a further subsequence,
\begin{equation*}
U_N\to U,\quad\text{in }H^1_{\mathrm{loc}}(\Omega)\cap C_{\mathrm{loc}}(\Omega),
\end{equation*}
and the norm convergence gives, for each $N$,
\begin{equation*}
\|U\|_{C(\overline{\Omega_N})}\leq C\epsilon^{9/8}.
\end{equation*}
Sending $N$ to infinity shows that $U\in C(\overline{\Omega})$ with
$\|U\|_{C(\overline{\Omega})}\leq C\epsilon^{9/8}$.

\noindent\textit{5}. Since $\nabla p_N$ is uniformly bounded by $C\epsilon^{5/4}$ in $L^2(\Omega_N)$, a diagonal argument yields, along a further subsequence,
\begin{equation*}
\nabla p_N\rightharpoonup\nabla p,\quad\text{in }L^2_{\mathrm{loc}}(\Omega),
\end{equation*}
for some $p\in H^1_{\mathrm{loc}}(\Omega)/\mathbb{R}$, and lower semicontinuity gives
\begin{equation*}
\|\nabla p\|_{L^2(\Omega_N)}\leq C\epsilon^{5/4}.
\end{equation*}
Sending $N$ to infinity shows that $\nabla p\in L^2(\Omega)$ with
$\|\nabla p\|_{L^2(\Omega)}\leq C\epsilon^{5/4}$.

\noindent\textit{6}. Let $\Phi\in C_c^\infty(\Omega;\mathbb{R}^2)$ and $\varphi\in C_c^\infty(\Omega;\mathbb{R})$. For sufficiently large $N$, both test functions $\Phi$ and $\varphi$ are supported in $\Omega_N$. The strong $H^1_{\mathrm{loc}}$ convergence of $U_N$ and weak $L^2_{\mathrm{loc}}$ convergence of $\nabla p_N$ justify passing to the limit in the weak formulation:
\begin{equation*}
\begin{cases}
\displaystyle\epsilon\int_\Omega\nabla U:\nabla\Phi+\int_\Omega((U^a+U)\cdot\nabla U)\cdot\Phi+\int_\Omega(U\cdot\nabla U^a)\cdot\Phi+\int_\Omega\nabla p\cdot\Phi=-\int_\Omega R^a\cdot\Phi,\\
\displaystyle\int_\Omega U\cdot\nabla\varphi=0.
\end{cases}
\end{equation*}
Since $U\in H^2(\Omega)$ and $\nabla p\in L^2(\Omega)$, the equations actually hold in the strong sense.

The boundary condition $U(x,0)=0$ is satisfied since $U\in H^1_0(\Omega)$. Moreover, $U$ is uniformly continuous because $H^2(\Omega)\hookrightarrow C_b^{0,\gamma}(\overline{\Omega})$ for $0<\gamma<1$, which, combined with $L^2$ integrability, implies the far-field condition $\lim\limits_{(x,y)\to\infty}U(x,y)=0$. This concludes the proof of Proposition \ref{global-nonlin-eq-exist}.
\end{proof}

\section{Proof of Main Theorem}
Finally, we conclude the main theorem.

\begin{proof}[Proof of Theorem \ref{main-thm}]
Theorem \ref{main-thm} is a direct consequence of the construction of the approximate solution \eqref{approx-profile}, Proposition \eqref{approx-profile-prop}, and Proposition \ref{global-nonlin-eq-exist}.
\end{proof}

\bigskip

\noindent\textbf{Acknowledgments}: A. L. Mazzucato was partially supported by the US National Science Foundation under grants DMS-2511023.
D. Wang was supported in part by the NSF grant DMS-2510532. Y. Wang was partially supported by the National Key Research and Development Program of China, grant 2021YFA1000800; the National Natural Science Foundation of China, grants 12288201 and 12421001; and the CAS Project for Young Scientists in Basic Research, grant YSBR-031.
\medskip

\noindent\textbf{Declaration of generative AI and AI-assisted technologies in the manuscript preparation process:}
During the preparation of this work, the authors used ChatGPT in order to perform minor language editing and reference checking. After using this tool/service, the authors reviewed and edited the content as needed and take full responsibility for the content of the published article.

\end{document}